\pdfoutput=1
\documentclass[a4paper,11pt]{article}
\usepackage{amsmath,amssymb,latexsym} 
\usepackage{times,cite,theorem}
\usepackage[]{graphicx}
\usepackage{ifthen}
\usepackage{rotating} 
\usepackage{array}
\usepackage[thinspace,thinqspace,squaren,textstyle]{SIunits}
\usepackage{color}

\def\ord{{\cal O}}\def\CO{\ord} \def\eps{\varepsilon} \def\rey{\text{R}}
\def\Rc{\rey_\text{crit}} \def\web{\text{W}}
\def\ibond{{\rm B}_{{\rm i}}}
\def\pair{p_\text{air}}
\def\Pair{P_\text{air}} \def\coloneqq{:=}

\def\bce{\begin{centering}}\def\ece{\end{centering}}
\newcommand{\reff}[1]{(\ref{#1})}

\def\al{\alpha}\def\del{\delta}\def\th{\theta}
\def\pa{\partial}
\newcommand{\bpm}{\begin{pmatrix}}\newcommand{\epm}{\end{pmatrix}}
\def\uhat{\hat{u}}\def\what{\hat{w}}
\def\res{\operatorname{Res}}
\newcommand{\spr}[1]{\left\langle #1 \right\rangle}

\def\Uprof{\tilde{U}} \def\Wprof{\tilde{W}} \def\QI{Q^{\text{IBL}}}
\def\bhat{\hat{b}} \def\Bhat{\hat{B}} \def\ahat{\hat{a}}
\def\lambdahat{\hat{\lambda}} \def\hhat{\hat{h}} 
\def\xhat{\hat{x}}\def\zhat{\hat{z}}
\def\Xhat{\hat{X}}\def\exhat{\mathbf{e}_{\hat{x}}}
\def\ezhat{\mathbf{e}_{\hat{z}}} \def\ex{\mathbf{e}_x} \def\ez{\mathbf{e}_z}
\def\ig{\includegraphics}

\newcommand{\dx}[1][]{\ifthenelse{\equal{#1}{}}{\partial_x}{\partial_x^{#1}}}
\newcommand{\dxhat}[1][]{\ifthenelse{\equal{#1}{}}{\partial_{\xhat}}{\partial_{\xhat}^{#1}}}
\newcommand{\dX}[1][]{\ifthenelse{\equal{#1}{}}{\partial_X\!}{\partial_X^{#1}\!}}
\newcommand{\dXhat}[1][]{\ifthenelse{\equal{#1}{}}{\partial_{\Xhat}\!}{\partial_{\Xhat}^{#1}\!}}
\newcommand{\dz}[1][]{\ifthenelse{\equal{#1}{}}{\partial_z}{\partial_z^{#1}}}
\newcommand{\dZ}[1][]{\ifthenelse{\equal{#1}{}}{\partial_Z\!}{\partial_Z^{#1}\!}}
\newcommand{\dt}[1][]{\ifthenelse{\equal{#1}{}}{\partial_t}{\partial_t^{#1}}}
\newcommand{\dT}[1][]{\ifthenelse{\equal{#1}{}}{\partial_T\!}{\partial_T^{#1}\!}}
\newcommand{\abt}[2]{\frac{d #1}{d #2}}

\begin{document}
\title{An integral boundary layer equation for film flow over inclined wavy
  bottoms}
 \author{T. H\"acker$^1$, H. Uecker$^2$\\[2mm]
\small $^1$Institut f\"ur Analysis, Dynamik und Modellierung, 
Universit\"at Stuttgart, Pfaffenwaldring 57,\\ 
\small D--70569 Stuttgart, tobias.haecker@mathematik.uni-stuttgart.de\\[2mm]
\small $^2$Institut f\"ur Mathematik, Carl von Ossietzky Universit\"at Oldenburg, \\
\small D--26111 Oldenburg, hannes.uecker@uni-oldenburg.de }
\maketitle
\normalsize
\begin{abstract}
  We study the flow of an incompressible liquid film down a wavy incline.
  Applying a Galerkin method with only one ansatz function to the
  Navier--Stokes equations we derive a second order weighted residual 
integral boundary layer equation, which in particular may be used to 
describe eddies in the troughs of the wavy bottom. 
We present numerical results which 
  show that our model is qualitatively and quantitatively accurate in wide
  ranges of parameters, and we use the model to study some new phenomena, 
for instance the occurrence of a short wave instability 
(at least in a phenomenological sense) for laminar flows 
which does not exist over flat bottom. 
\end{abstract}
\section{Introduction}
The gravity driven free surface flow of a viscous incompressible fluid down an
inclined plate has various
engineering applications, for instance in cooling and coating processes. 
For a flat bottom the problem,  governed by the Navier--Stokes equations, 
is extensively studied experimentally, 
numerically and analytically, see, e.g., \cite{cd02} for a review. 
In particular it is well known that there exists a stationary solution
with a parabolic velocity profile and a flat surface, the so called Nusselt
solution, which is unstable to long waves if the Reynolds number
exceeds a critical value $\Rc = 5/6 \cot \alpha$, where $\alpha$ is
the inclination angle \cite{Benjamin_57, Yih_63}. 
However, the Navier--Stokes equations in combination with the free
surface are hard to handle and one is often not interested in the flow field
but only in, e.g., the film thickness $F$. Thus there has been much effort to
derive model equations for the evolution of $F$. Because of the long wave
character of the instability, length scales of free surface perturbations are
large compared to the film thickness. Therefore a small parameter $\eps$ can
be introduced to scale downstream derivatives. By an asymptotic expansion
approach a scalar evolution equation for $F$ was derived in \cite{Benney_66}
and later corrected in \cite{Lin_74}. However, this so called Benney equation
has finite-time blow-up solutions even at moderate Reynolds numbers, see
\cite{Pumir_83}. Nevertheless, asymptotically it can be used to check the
consistency of improved models, see \cite{Scheid_06}. 

Besides the reduction of the Navier--Stokes problem to a scalar equation 
for the film thickness $F$ a hierarchy of less drastic reductions 
has been studied, starting with so called boundary layer equations, 
see again \cite[Chapter 2]{cd02}, for instance. An important step was the derivation of an integral boundary layer equation (IBL) by Shkadov in
\cite{Shkadov_67}. He used the averaging method of K\'{a}rm\'{a}n--Pohlhausen
which consists of taking a parabolic velocity profile like the stationary
Nusselt solution as ansatz for the downstream velocity component $U$ and
integrating the streamwise momentum equation along the $Z$ coordinate
perpendicular to the bottom.  This yields a system of two evolution equations
for $F$ and the local flow rate $Q = \int_0^F U dZ$.

Although the IBL reproduces various experimental observations like the
existence of solitary waves it shows the following inaccuracies:
\begin{enumerate}
\item The predicted critical Reynolds number differs from the exact value by a
  factor $5/6$.
\item The IBL is not consistent with the Benney equation.
\item The assumed parabolic velocity profile does not fulfill the dynamic
  boundary condition at second order.
\end{enumerate}
The first problem follows from a linear stability analysis which
yields $\text{R}_{\text{crit, IBL}} = \cot \alpha$. For the second point one
derives a scalar evolution equation for $F$ from the IBL. This can be done by
enslaving the flow rate $Q$ to the film thickness $F$ and expanding it in
powers of $\eps$, which gives a scalar equation for $\dT F$ differing from
the Benney equation already at order $\eps$, see \cite{Ruyer_98}. The third
problem is due to the fact that the parabolic velocity profile has its maximum
at the free surface which implies $\dZ U(F) = 0$.

Recently there has been much effort to overcome these
problems. Along \cite{Ruyer_98, Ruyer_00, Scheid_06} a
two-equation model for $F$ and $Q$ has been derived by a Galerkin method. 
Based again on a long wave expansion of the Navier--Stokes
equations, the Nusselt solution and three more 
polynomials appearing in the derivation of the Benney
equation served as ansatz and test functions. 
The resulting model consisted of four
evolution equations for $F, Q$ and two other quantities measuring the
deviation from the parabolic velocity profile. From this a simplified model, 
called weighted residual integral boundary layer equation (WRIBL) 
for $F$ and $Q$ was derived which is consistent with the
Benney equation at order $\eps^2$ and predicts the correct critical Reynolds
number. However, this model does not reproduce well known solitary wave solutions if the Reynolds number exceeded a certain value only slightly larger than the instability threshold. This deficiency can be cured by a Pad\'{e}-like regularization method in \cite{Scheid_06}. Moreover, in numerical simulations the extension of the WRIBL to three-dimensional flows yields excellent agreement with recent experimental results 
from \cite{pano03}, see again \cite{Scheid_06}. See also \cite{ogn08} for further detailed numerical studies of this model.

The problem over wavy bottom is studied much less extensively. For
experimental results we refer to 
\cite{poz88,vb02,wsa03,Acta_05, vmhm05, avb06, wbhua08}. 
On the theoretical side, \cite{wa03,Acta_05} give an expansion of Nusselt 
like stationary solutions in suitable small parameters and an analysis 
of their stability. In \cite{Trifonov_98,trif04,trif07a} the problem 
is studied numerically by simulations 
of both the full Navier--Stokes problem and model equations 
  derived in a similar way as in \cite{Shkadov_67}. 
Moreover, a detailed numerical stability analysis based on the 
Navier--Stokes equations has been carried out \cite{trif07}.  
In \cite{dav07} a scalar Benney like model has been derived and 
studied numerically, and 
in \cite{hbaw08} an IBL over wavy bottom has been derived using 
Shkadov's method.  Finally, using the method 
from \cite{Ruyer_00,Scheid_06} 
a first-order WRIBL has been derived and studied in great detail 
in \cite{oh08}. 
 
Here we continue into a similar direction as \cite{oh08} 
by deriving and analyzing 
numerically an alternative WRIBL equation and a regularized version. 
However, in contrast to \cite{oh08} our analysis is based on 
curvilinear coordinates 
from \cite{wsa03} which allow to treat more general 
situations where for instance the free surface is not necessarily a graph 
over the (flat bottom) downstream coordinate. These curvilinear coordinates 
are also more natural since they allow a clear distinction between 
flow components tangential and normal to the bottom. 
Moreover, 
our  WRIBL is second order accurate which for instance 
allows the description of eddies in the troughs of the wavy bottom. 
Finally, our approach is somewhat simpler than the (more general) 
approach of \cite{Ruyer_00,Scheid_06} which consists of several 
polynomial ansatz and test functions in the Galerkin expansion.  
We find that by taking an accurate velocity profile $\Uprof$
as single ansatz and test function in the Galerkin method the WRIBL  
can be obtained in one step. 

Thus, the outline is as follows: In Section 2 we present the governing equations
in curvilinear coordinates. Since we focus on film flow over bottoms with 
long wave undulations we assume the bottom steepness and the 
non-dimensional wave number to be of order $\eps$, $0<\eps\ll 1$, and expand 
all equations up to $\ord(\eps^2)$. In
Section 3 we derive an appropriate velocity profile serving as ansatz and test
function used to derive our WRIBL by the Galerkin method in Section 4, and 
in Section 5 we check the consistency of the resulting WRIBL with the 
Benney equation over wavy bottoms. From the WRIBL we derive a regularized 
version called rWRIBL in Section 6 by removing second-order 
inertia terms which otherwise may lead to some unphysical behaviour.
In Section 7 we finally give some numerical results. First, in \S\ref{onum}, 
by comparison with available 
experimental and full Navier--Stokes numerical data  we illustrate the 
accuracy of our rWRIBL over wide parameter regimes, including the occurrence 
of eddies. Second, in \S\ref{nnum} we illustrate 
two new phenomena, namely that the bottom modulation may introduce a 
short wave instability (in a phenomenological sense) not present over 
flat bottom (except for rather 
extreme parameter ranges), and that and how the free surface may cease to be 
a graph over the (flat bottom) downstream coordinate. A short summary 
is given in \S\ref{ssec}.  
\section{Governing equations}
\begin{figure}
\begin{center}
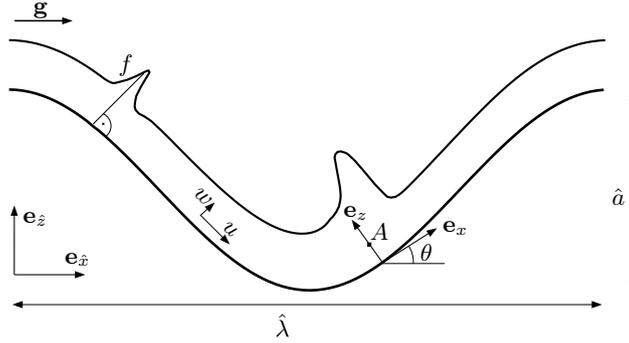
\caption{Sketch of the geometry and the curvilinear coordinate system.}
  \label{fig1}
\end{center}
\end{figure}

Figure \ref{fig1} illustrates the inclined film problem 
with an undulated bottom $\bhat$.  The liquid is assumed incompressible and
Newtonian, the Cartesian coordinate system $\exhat, \ezhat$ is
inclined at an angle $\alpha$ with respect to the horizontal ($\al=90^\circ$ 
in Fig.~\ref{fig1}), and the bottom profile
$\bhat(\xhat)$ is periodic with wavelength $\lambdahat$ and amplitude
$\ahat$. 
As we want to expand the governing equations in a small parameter 
$\eps$ it is useful and natural to introduce a
curvilinear coordinate system for the following reasons. 
First, although
the Nusselt solution is no longer a stationary solution if the bottom is
undulated, for thin films and low Reynolds numbers the flow 
 $(u,w)$ is still mainly parallel to the bottom. 
To apply different scalings to $u$ and $w$ the coordinate
system thus has to be orientated along the bottom profile such that the
$u$ component is tangential to the bottom, while using 
a fixed Cartesian coordinate system scaling involves a mixing of the 
Cartesian velocity components $\uhat$, $\what$. 
 Second, for larger 
Reynolds numbers we may anticipate situations as sketched in 
Fig.~\ref{fig1} where the free surface is not a graph over 
$\xhat$ and cannot easily be described in Cartesian coordinates. 

Thus, at every point of the bottom $\xhat
\exhat+\bhat(\xhat)\ezhat$ we define a local coordinate system $\ex, \ez$ with
$\ex$ tangential and $\ez$ normal to the bottom. 
For an arbitrary point $A$
within the liquid the arc length $x$ of the bottom 
 and the distance $z$ along $\ez$ to the bottom are now 
taken as curvilinear coordinates.  
As we focus on film flow over weakly undulated bottoms this relation is always unique. Thus, 
$$
A=\bpm \xhat-\sin\th\, z\\ \bhat(\xhat)+\cos\th\, z\epm
$$
in $\exhat$, $\ezhat$ coordinates, 
where $\th=\th(x)$ is the local inclination angle between $\exhat$ and $\ex$. 
In order to transform gradients we will also need the
bottom curvature $\kappa$ which is defined by
\begin{equation}
  \kappa(\xhat) = -\frac{\dxhat[2] \bhat(\xhat)}
{(1 + (\dxhat \bhat(\xhat))^2)^{\frac{3}{2}}}.
  \label{kappa}
\end{equation}
For further details concerning the transformation to curvilinear coordinates
we refer to \cite{Acta_05}.

To describe the free-surface flow we introduce the variables in 
Table \ref{tab1}. 
\begin{table}
\begin{center}
  \begin{tabular}{ll}
    $\mathbf{v}(x,z,t) = u(x,z,t) \ex+w(x,z,t) \ez$ & velocity field \\
    $f(x,t)$  & film thickness (perpendicular to the bottom) \\
    $p(x,z,t)$ & pressure \\
    $\pair$ & pressure of the air above the liquid surface \\
    $\sigma$ & surface tension \\
    $\rho$ & liquid density \\
    $ \nu$ & kinematic viscosity \\
    $\mathbf{g}$ & gravity acceleration
  \end{tabular}
\caption{Physical quantities.}\label{tab1}
\end{center}
\end{table}
In contrast to Cartesian coordinates all quantities measured in
curvilinear coordinates are written without a hat. The governing two-dimensional Navier--Stokes equations now read
\begin{align}
  \dt u&+\frac{1}{1+\kappa z}u \dx u+w \dz u+\frac{1}{1+\kappa z} \kappa uw
  \nonumber \\
  & =-\frac{1}{\rho} \frac{1}{1+\kappa z} \dx p+g \sin(\alpha-\theta)
  +\nu \left[ \frac{1}{(1+\kappa z)^3}\dx \kappa(w -z \dx u) \right. \nonumber \\
  &\qquad+\frac{1}{(1+\kappa z)^2}(\dx[2] u-\kappa^2 u+2 \kappa \dx w) +
  \left. \frac{1}{1+\kappa z} \kappa \dz u+\dz[2] u \right],
  \label{eq1}\\
  \dt w&+\frac{1}{1+\kappa z}u \dx w+w \dz w-\frac{1}{1+\kappa z} \kappa u^2 \nonumber \\
  & =-\frac{1}{\rho} \dz p-g \cos(\alpha-\theta)+\nu \left[ -\frac{1}{(1+\kappa z)^3}\dx \kappa(u +z \dx w) \right. \nonumber \\
  &\qquad+\frac{1}{(1+\kappa z)^2}(\dx[2] w-\kappa^2 w-2 \kappa \dx u) +
  \left. \frac{1}{1+\kappa z} \kappa \dz w+\dz[2] w \right],
  \label{eq2}\\
  \frac{1}{1+\kappa z} &(\dx u+\kappa w)+\dz w = 0.
  \label{eq3}
\end{align}
At the bottom $z \equiv 0$ we have the no-slip and no-flux condition
\begin{equation}
  u\big|_{z=0} = w\big|_{z=0} = 0.
  \label{eq4}
\end{equation}
The dynamic boundary condition tangential and normal to the free surface
$z\equiv f$ reads
\begin{align}
  0&= \left( (1+ \kappa f)^2-(\dx f) ^2 \right) \left( \frac{\pa_x w-\kappa u}{1+\kappa f}+\pa_z
  u \right)+4(1+ \kappa f)\dx f \dz w,
  \label{eq5}\\
  \sigma& \frac{(1+\kappa f) \dx[2]f-f \dx \kappa \dx f-\left((1+ \kappa
      f)^2+2(\dx f)^2 \right) \kappa}{\left( (1+ \kappa f)^2+(\dx f)^2
    \right)^{3/2}}+(p-\pair)\notag\\
&\qquad = \frac{2\rho\nu}{1+(\pa_x f/(1+\kappa f))^2}
\left(\frac{(\pa_x f)^2(\pa_x u+\kappa w)}{(1+\kappa f)^3}+\pa_z w
-\frac{\pa_x f}{1+\kappa f}\left(\frac{\pa_x w-\kappa u}{1+\kappa f}+\pa_z
  u\right)
\right)
  \label{eq6}
\end{align}
while the kinematic boundary condition is
\begin{equation}
  \abt{}{t}(f(x,t)-z) = 0 \quad \Leftrightarrow \quad \dt f+\frac{1}{1+\kappa f} u \dx f-w = 0.
  \label{eq7}
\end{equation}
In order to introduce dimensionless quantities we refer to the 
stationary solution over a flat incline. 
This so called Nusselt solution has the mean flow velocity
$\langle u \rangle = \frac{g \sin \alpha \hhat^2}{3 \nu}$, where $\hhat$ is
the constant film thickness.  We set 
\begin{alignat*}{4}
  & X = \frac{2\pi}{\lambdahat} x, & & Z = \frac{1}{\hhat}z, & & F = \frac{1}{\hhat} f, & \quad & U = \frac{1}{\langle u \rangle}u,\\
  & W = \frac{\lambdahat}{2 \pi \hhat\langle u \rangle}w, & \quad & T =
  \frac{2 \pi \langle u\rangle}{\lambdahat}t, &\quad & K =
  \frac{\lambdahat^2}{4\pi^2 \ahat}\kappa, & \quad& P = \frac{1}{\rho \langle
    u \rangle^2}p.
\end{alignat*}
Additional to $\alpha$ we can choose four non-dimensional parameters to write
the governing equations dimensionless. To describe surface tension and
viscosity effects we use
\begin{align*}
 \ibond&\coloneqq \frac{4\pi^2 l^2_{\text{ca}}}{\lambdahat^2
    \sin\alpha}=\frac{4 \pi^2\sigma}{\rho g \lambdahat^2 \sin \alpha}
  \quad\text{(inverse Bond number), } \\
 \rey&\coloneqq \frac {\langle
    u\rangle\hhat}{\nu} =\frac{g \hhat^3 \sin \alpha}{3 \nu^2}\quad
  \text{(Reynolds number).}
\end{align*}
Here $l_{\text{ca}} = \left( \frac{\sigma}{\rho g} \right)^{\frac{1}{2}}$ is
the capillary length. The relation of $\ibond$ to the also frequently used
Weber number $\web=\frac{\sigma}{\rho g \hhat^2 \sin\alpha}$ is
$\web=\frac{1}{\delta^2}\ibond.$ For the geometric quantities we introduce
\begin{align*}\delta := 2 \pi \frac{\hhat}{\lambdahat}\quad
  \text{(dimensionless wave number),}\qquad \zeta :=2 \pi
  \frac{\ahat}{\lambdahat}\quad\text{(bottom steepness).}
\end{align*}
As we are interested in thin films over weakly undulated bottoms we suppose
throughout that both $\delta$ and $\zeta$ are of order $\eps$, where $\eps$ is
a small parameter, while $\rey, \ibond$ and $\alpha$ are
assumed to be of order $1$. The latter means that $\alpha$ is bounded 
away from zero such that $\cot(\al)$ is bounded. However, $\al=90^\circ$ 
such that $\cot(\al)=0$ is allowed. 

{\em All calculations will be exact of order $\eps^2$},
i.e.~we keep all terms of order $1, \delta, \zeta, \delta^2, \zeta^2$ and
$\delta \zeta$. Throughout we will only display the $\ord(\eps^3)$-symbol if
we want to emphasize that our calculations are only asymptotically correct. In
all other cases we will skip it.  In particular, skipping $\CO(\eps^3)$-terms,
the dimensionless governing equations read
\begin{align}
  & \delta\rey \dT U+\delta\rey \dX U U+\delta\rey \dZ U W = -\delta\rey \dX P
  + 3\frac{\sin(\alpha{-}\theta)}{\sin\alpha} +\delta^2 \dX[2] U+\delta\zeta K
  \dZ U+\dZ[2] U,
  \label{NS1}\\
  & \delta^2\rey \dT W+\delta^2\rey U \dX W+\delta^2 \rey \dZ W W -
  \delta\zeta\rey K U^2 = -\rey \dZ P - 3\frac{\cos(\alpha{-}\theta)}
{\sin\alpha} +\delta \dZ[2] W,
  \label{NS2}\\
  &\dX U+\dZ\left( (1+\delta\zeta KZ)W \right) = 0, 
  \label{NS3}\\
  &U(0) = W(0) = 0, 
  \label{NS4}\\
  & (1{+}2\delta\zeta KF{-}\delta^2 (\dX F)^2) \dZ U(F)+\delta^2 \dX W(F)-\delta
  \zeta K U(F)+4 \delta^2 \dX F \dZ W(F) = 0,
  \label{NS5}\\
  & 3 \ibond (\dX[2] F-\xi K) = -\rey(P(F)-\Pair)+2 \delta\dZ W(F)+\CO(\eps^2),
  \label{NS6}\\
  &\dT F+(1-\delta\zeta KF)\dX F U(F)-W(F) = 0.
  \label{NS7}
\end{align}
The dynamic boundary condition normal to the free surface \eqref{NS6}, where
we used the abbreviation $\xi\coloneqq\frac{\zeta}{\delta}$, is only given up
to 
order $\eps$. As we are not interested in second-order terms of the pressure
$P$ this turns out to be sufficient.
\section{A first-order velocity profile}
\label{sec_prof}
For given $F$ we derive a solution $(U,W,P)$ of the time dependent equations
\eqref{NS1}--\eqref{NS6} which is exact to order $\eps$. By introducing the
flow rate $Q$ as independent quantity we also construct a velocity
profile $\Uprof$ which will serve as ansatz and test function in the Galerkin
approach in Section \ref{sec_gal}. There, a first-order profile
$\Uprof=\Uprof_0 + \eps \Uprof_1$ is sufficient since we can extract all
necessary second-order terms from the boundary conditions.

We assume that $F$ is of order $1$ while the velocity field $(U,W)$ and the
pressure $P$ are enslaved by $F$ and can be expanded in powers of $\eps$:
\begin{equation}
  U = U_0+\eps U_1+\ord(\eps^2),\quad 
  W = W_0+\eps W_1+\ord(\eps^2),\quad 
  P = P_0+\eps P_1+\ord(\eps^2).
  \label{3_1}
\end{equation}
The geometric quantities $K$ and $\theta$ coming from the bottom profile can
be expanded in powers of $\eps$, too. It turns out that the bottom curvature
$K$ does not contain terms of first order while the local inclination angle
has a leading $\zeta$, i.e.
\begin{equation*}
  K = K_0 + \zeta^2 K_2 + \ord(\zeta^4), \quad \theta = 
\zeta\theta_1 + \ord(\zeta^3)
\end{equation*}
with $\theta_1(X)=\pa_X \hat{B}(X)$, see Appendix \ref{sec_CC}. This yields
\begin{equation*}
  \frac{\cos(\alpha{-}\theta)}{\sin\alpha} 
= \cot\alpha + \zeta\theta_1 - \frac{1}{2}\zeta^2\cot\alpha 
\theta_1^2 + \ord(\zeta^3), \quad 
\frac{\sin(\alpha{-}\theta)}{\sin\alpha} 
=1-\zeta\cot\alpha \ \theta_1 - \frac{1}{2}\zeta^2\theta_1^2 + \ord(\zeta^3).
\end{equation*}
Since both $\delta$ and $\zeta$ are of order $\eps$, equations
\eqref{NS1}--\eqref{NS6} read at $\ord(1)$
\begin{align*}
  &3+\dZ[2] U_0 = 0,\quad -\rey \dZ P_0-3 \cot \alpha = 0, \quad
  \dX U_0+\dZ W_0 = 0, \\
  &U_0(0) = W_0(0) = 0, \quad \dZ U_0(F) = 0, \quad 3 \ibond(\dX[2] F-\xi K_0)
  = -\rey(P_0(F)-\Pair).
\end{align*}
The $\ord(1)$-solution thus is
\begin{align}
  &U_0 = -\frac{3}{2}Z^2+3FZ,\label{3_17}
  \quad W_0 = -\frac{3}{2}\dX F Z^2, \quad P_0 = \frac{3}{\rey}(\cot \alpha
  (F-Z)-\ibond \dX[2]F+\ibond \xi K_0) +\Pair.
\end{align}
At $\ord(\eps)$ we get the equations
\begin{align*}
  &\delta\rey \dT U_0+\delta\rey\dX U_0 U_0+\delta\rey\dZ U_0 W_0 = -\delta\rey\dX P_0 - 3\zeta\cot\alpha \ \theta_1 + \eps\dZ[2]U_1, \\
  &-\eps\rey\dZ P_1 -3\zeta\theta_1 +\delta\dZ[2]W_0 = 0, \\
  &\dX U_1+\dZ W_1 = 0, \quad 
  U_1(0) = W_1(0) = 0, \quad \dZ U_1(F) = 0, \quad -\eps\rey
  P_1(F)+2\delta\dZ W_0(F) = 0, \notag
\end{align*}
with solutions
\begin{alignat}{2}
  &\eps U_1 =&\ &\frac{1}{2}\delta\rey \dT F(Z^3{-}3F^2 Z)+\delta\rey\dX 
F\left(\frac{3}{8}FZ^4-\frac{3}{2}F^4 Z\right)  \nonumber \\
  &    &   & + 3(\delta\cot\alpha \dX F - \delta \ibond\dX[3]F 
+ \zeta\ibond\dX K_0 + \zeta\cot\alpha \theta_1) 
\left( \frac{1}{2}Z^2{-}FZ \right), \label{3_19}\\
  &\eps W_1 =&\ &-\frac{1}{2}\delta\rey\pa_{TX}\!F
\left( \frac{1}{4}Z^4-\frac{3}{2}F^2Z^2 \right)+\frac{3}{2}\delta\rey
\dT F \dX F F Z^2\nonumber \\
&&&-\delta\rey \dX[2]F \left( \frac{3}{40}FZ^5-\frac{3}{4}F^4Z^2 \right) 
  -\delta\rey(\dX F)^2 \left( \frac{3}{40}Z^5-3F^3Z^2 \right) 
\notag\\ 
&&&+ \frac{3}{2}(\delta\cot\alpha \dX F - \delta \ibond \dX[3]F 
+ \zeta\ibond\dX K_0 + \zeta\cot\alpha \ \theta_1)\dX F Z^2 \nonumber \\
  &    &   & - (\delta \cot\alpha \dX[2]F - \delta\ibond\dX[4]F 
+ \zeta \ibond \dX[2] K_0 + \zeta\cot\alpha\dX \theta_1)\left( \frac{1}{2}Z^3 
- \frac{3}{2}FZ^2 \right), \nonumber \\
  &\eps P_1 =& &-\frac{3}{\rey}\zeta\theta_1(Z-F) - \frac{3}{\rey}\delta\dX
  F(Z+F). \nonumber
\end{alignat}
To get rid of the time derivatives of $F$ we use the kinematic boundary
condition \eqref{NS7} which leads at $\ord(1)$ to the identity
\begin{align}
  \dT F & = -\dX F U_0(F)+W_0(F)+\ord(\eps) = -3 \dX F F^2+\ord(\eps).\notag
\end{align}
Thus $U_1$ can be rewritten as
\begin{align}
  \eps U_1 = & \delta\rey \dX F F^5 \left(
    \frac{3}{8}\left(\frac{Z}{F}\right)^4
-\frac{3}{2}\left(\frac{Z}{F}\right)^3+3\frac{Z}{F} \right)\nonumber \\
  & -3 F^2(\delta \cot\alpha \dX F{-}\delta\ibond\dX[3]
F{+}\zeta\ibond\dX K_0{+}\zeta\cot\alpha \ \theta_1) 
\left(\frac{Z}{F}{-}\frac{1}{2}\left( \frac{Z}{F}\right)^2 \right).
  \label{3_18}
\end{align}

If we assume temporarily that also the local flow rate $Q = \int_0^F UdZ$ is
enslaved by $F$ we can easily state the $\eps$-expansion of $Q = Q_0+\eps Q_1
+ \ord(\eps^2)$, namely
\begin{align}
  Q_0=&\int_0^F U_0 dZ = F^3, \label{3_28} \\
  \eps Q_1 =& \eps\int_0^F U_1 dZ=
\frac{6}{5} \delta\rey \dX F F^6 -
  F^3(\delta\cot\alpha \dX F{-}\delta\ibond\dX[3]F{+}\zeta\ibond\dX K_0{+}
  \zeta\cot\alpha \theta_1).
  \label{3_29}
\end{align}
As mentioned in the introduction we cannot maintain the enslavement of $Q$ to
$F$ since this would lead to a single evolution equation for $F$ which fails
to reproduce physics correctly. Therefore we treat $Q$ as independent
$\ord(1)$-quantity and introduce a second representation
\begin{gather}
  \Uprof =\Uprof(F,Q)= \Uprof_0+\eps\Uprof_1+\ord(\eps^2)\label{2nd}
\end{gather}
of the velocity profile which depends on both $F$ and $Q$. For consistency, if
we plug the enslaved version $Q = Q_0+\eps Q_1+\ord(\eps^2)$ into \reff{2nd}
we must recover the expansion $U = U_0+\eps U_1+\ord(\eps^2)$ calculated in
\eqref{3_17}, \eqref{3_18}. This yields the following conditions for
$\Uprof_0$:
\begin{gather}
  \text{(i)}\quad \int_0^F \Uprof_0 dZ = Q \quad\text{as } Q\text{ is of order
  }1,\qquad\text{(ii)}\quad \Uprof_0 = U_0\quad\text{if } Q = Q_0 +
  \ord(\eps)\text{ is assumed.}  \notag\end{gather} As $Q$ is independent of
$Z$ the first condition implies that $Q$ occurs as a factor in
$\Uprof_0$. From \eqref{3_28} we know that in the enslaved version of $Q$ in
zeroth order we have $Q=F^3$.  Thus
\begin{equation}
  \Uprof_0 = \frac{3Q}{F}\left( -\frac{1}{2} 
    \left(\frac{Z}{F}\right)^2+\frac{Z}{F} \right), 
\end{equation}
which is exactly the lubrication ansatz which is used in the method of
K\'{a}rm\'{a}n--Pohlhausen. Thus our new velocity profile will emerge as
refinement of the parabolic profile.

On the other hand, plugging $Q = Q_0+\eps Q_1$ into $\Uprof_0$ yields
\begin{align}
  \Uprof_0 = \ & -\frac{3}{2}Z^2+3FZ +\delta\rey \dX F F^5 \left( -\frac{9}{5}\left(\frac{Z}{F}\right)^2+\frac{18}{5}\frac{Z}{F} \right) \nonumber\\
  &-3F^2(\delta\cot\alpha \dX F - \delta\ibond\dX[3]F + \zeta\ibond\dX K_0 +
  \zeta\cot\alpha \
  \theta_1)\left(-\frac{1}{2}\left(\frac{Z}{F}\right)^2+\frac{Z}{F}
  \right). \label{3_29b}
\end{align}
Thus, comparing \eqref{3_18} and \eqref{3_29b}, $\Uprof_0$ contains terms which belong to $U_1$, and therefore $\Uprof_1$
consists of less terms than $U_1$, namely
\begin{equation}
  \eps\Uprof_1 = \delta\rey \dX Q Q \left( \frac{1}{8}
    \left(\frac{Z}{F}\right)^4-\frac{1}{2}\left(\frac{Z}{F}\right)^3
    +\frac{3}{5}\left(\frac{Z}{F}\right)^2-\frac{1}{5}\frac{Z}{F} \right).
  \label{uti1}
\end{equation}
To sum up, if $Q$ is treated as independent $\ord(1)$-quantity we obtain the
first-order velocity profile
\begin{equation}
  \Uprof = \frac{3Q}{F}\left(\frac{Z}{F}{-}\frac{1}{2}
    \left(\frac{Z}{F}\right)^2\right)
+\delta\rey \dX Q Q \left( \frac{1}{8}\left(\frac{Z}{F}\right)^4
{-}\frac{1}{2}\left(\frac{Z}{F}\right)^3{+}\frac{3}{5}\left(\frac{Z}{F}\right)^2
{-}\frac{1}{5}\frac{Z}{F} \right).
  \label{3_40}
\end{equation}
Similarly, the second-order velocity profiles $U_2$ and $\Uprof_2$ are derived
in Appendix \ref{sec_u2}. These are not needed for the derivation of the WRIBL
but for the reconstruction of the flow field in Section
\ref{sec_num}.
\section{Galerkin method}
\label{sec_gal}

We start with the derivation of the evolution equation for $F$ by integrating
the continuity equation \eqref{NS3} along $Z$, i.e.
\begin{equation}
  \int_0^F \dX U dZ+\left[(1+\delta\zeta KZ)W\right]_0^F = 0.\notag
\end{equation}
From $Q = \int_0^F U dZ$ and the no-flux boundary condition we obtain $\dX
Q-\dX F U(F)+(1+\delta \zeta KF)W(F)=0$, and eliminating $W(F)$ by the
kinematic boundary condition \eqref{NS7} and skipping all terms of order
$\eps^3$ and higher finally gives
\begin{equation}
  \dT F = -(1-\delta\zeta KF)\dX Q. \label{FT}
\end{equation}
In order to derive an evolution equation for $Q$ we first eliminate the
pressure $P$ from the streamwise momentum equation \eqref{NS1} before we apply
a Galerkin method. By means of \eqref{NS2} $P$ can be written as
\begin{align}
  \delta\rey P(Z) & = \delta\rey P(F)-\delta\rey\int_Z^F \dZ P dZ \nonumber \\
  & = \delta\rey
  P(F)+3\delta\frac{\cos(\alpha-\theta)}{\sin\alpha}(F-Z)-\delta^2(\dZ
  W(F)-\dZ W(Z)).\notag
\end{align}
To eliminate $P(F)$ we use the dynamic boundary condition normal to the free
surface \eqref{NS6} and the continuity equation \eqref{NS3} to obtain
\begin{align}
  \delta\rey P(Z) & =  \delta\rey \Pair+\delta^2(\dZ W(F)+\dZ W(Z))
{-}3\ibond(\delta \dX[2]F{-}\zeta K)+3\delta\frac{\cos(\alpha{-}\theta)}{\sin\alpha} (F{-}Z)\notag\\
  & = \delta\rey\Pair{-}\delta^2(\dX U(F)+\dX U(Z)){-}3\ibond(\delta \dX[2]F{-}\zeta
  K)+3 \delta\frac{\cos(\alpha{-}\theta)}{\sin\alpha}(F{-}Z). \notag
\end{align}
Plugging this into the streamwise momentum equation \eqref{NS1} we obtain
\begin{align}
  & \delta\rey \dT U+\delta\rey \dX U U+\delta\rey W \dZ U \nonumber \\
  & = 3\frac{\sin(\alpha{-}\theta)}{\sin\alpha}+\dZ[2]U+2\delta^2 \dX[2]U
+3\delta \ibond \dX[3]F-3\zeta\ibond\dX K
-3\delta\frac{\cos(\alpha{-}\theta)}{\sin\alpha}\dX F 
\nonumber \\
  &\quad-3\delta\frac{\sin(\alpha{-}\theta)}{\sin\alpha}\dX\theta(F{-}Z) +
  \delta^2\frac{d}{dX}(\dX U(F))+\delta\zeta K \dZ U. \label{4_10}
\end{align}

The next step is to perform a Galerkin method with the single test and ansatz
function $\Uprof$ from \eqref{3_40}. Thus we plug $\Uprof$ into \eqref{4_10},
multiply the residual by $\Uprof$ itself and integrate the result along
$Z$. {\em We want all calculations to be exact of order $\eps^2$}. 
This seems to be
a problem since the first two terms on the right-hand side of \eqref{4_10} are
of order $1$ and we know $\Uprof = \Uprof_0+\eps\Uprof_1 + \eps^2\Uprof_2$
only up to $\ord(\eps)$. However, the first term
$3\frac{\sin(\alpha-\theta)}{\sin\alpha}$ is independent of $Z$, and by the
definition of $Q$ we get
\begin{equation}
  \int_0^F 3\frac{\sin(\alpha-\theta)}{\sin\alpha} \Uprof dZ =
  3\frac{\sin(\alpha-\theta)}{\sin\alpha}Q.\notag
\end{equation}
The second term $\dZ[2]U$ is slightly harder to manage. Integration by parts
together with the no-slip condition \mbox{$\Uprof(0) = 0$} yields 
\begin{equation}
  \int_0^F\dZ[2]\Uprof \Uprof dZ = \dZ\Uprof(F) \Uprof(F)-\int_0^F(\dZ \Uprof)^2dZ, \label{4_11}
\end{equation}
and up to order $\eps^2$ the integral on the right-hand side reads
\begin{align}
  \int_0^F(\dZ \Uprof)^2 dZ & = \int_0^F\left( (\dZ\Uprof_0)^2+2\eps\dZ
    \Uprof_0 \dZ \Uprof_1
+\eps^2(\dZ\Uprof_1)^2+2\eps^2\dZ\Uprof_0\dZ\Uprof_2 \right) dZ \nonumber \\
  & = 3\frac{Q^2}{F^3}+\frac{1}{175} \delta^2\rey^2\frac{1}{F}(\dX Q)^2 Q^2
  +6\eps^2\frac{Q}{F}\int_0^F\left( \frac{1}{F}-\frac{Z}{F^2} \right)
  \dZ\Uprof_2 dZ. \label{4_21}
\end{align}
At this point we need some information about the second-order term
$\eps^2\Uprof_2$. The velocity profile $\Uprof$ emanates from the asymptotic
solution $U$, and thus fulfills the boundary conditions
\eqref{NS4}, \eqref{NS5}. Moreover, $\int_0^F \Uprof_0 dZ = Q$, which implies
$\int_0^F \Uprof_2 dZ = 0$. Therefore and due to the no-slip boundary
condition the last integral in \eqref{4_21} satisfies
\begin{equation}
  \int_0^F\left( \frac{1}{F}-\frac{Z}{F^2} \right) \dZ\Uprof_2 dZ 
  = \left[ \left( \frac{1}{F}-\frac{Z}{F^2} \right) \Uprof_2 \right]_0^F
  +\frac{1}{F^2}\int_0^F \Uprof_2 dZ = 0\notag
\end{equation}
which gives
\begin{equation}
  \int_0^F(\dZ \Uprof)^2 dZ = 3\frac{Q^2}{F^3}+\frac{1}{175} 
  \delta^2\rey^2 \frac{1}{F} (\dX Q)^2 Q^2.\notag
\end{equation}
It remains to calculate the first term on the right-hand side of
\eqref{4_11}. From 
\eqref{NS5} we know that $\dZ \Uprof(F) =-\delta^2 \dX \Wprof(F) -
4\delta^2 \dX F \dZ \Wprof(F)+\delta\zeta K \Uprof(F)$ is of order $\eps^2$
where the velocity component $\Wprof$ can be expressed by $\Uprof$ due to the
continuity equation \eqref{NS3}. Thus the $\ord(1)$-terms of $\Uprof$ are
sufficient which means that we do not have to know $\Uprof_2$ explicitly. This leads finally to
\begin{align}
  \int_0^F\dZ[2] \Uprof \Uprof dZ = \ & \frac{3}{2} \delta^2\frac{1}{F}\dX[2]Q
  Q-\frac{9}{2}\delta^2\frac{1}{F^3}Q^2(\dX
  F)^2-\frac{9}{4}\delta^2\frac{1}{F^2}Q^2\dX[2]F
  +\frac{9}{2}\delta^2\frac{1}{F^2}\dX
  Q Q \dX F \nonumber \\
  &+\frac{9}{4}\delta\zeta K \frac{Q^2}{F^2}-3\frac{Q^2}{F^3} -
  \frac{1}{175}\delta^2\rey^2\frac{1}{F}(\dX Q)^2 Q^2.\notag
\end{align}
The other terms in \eqref{4_10} are all at least of order $\eps$ and we can
calculate them rather easily by plugging in $\Uprof = \Uprof_0 +
\eps\Uprof_1$. Testing \eqref{4_10} with $\Uprof$ leads to
\begin{align}
  \delta\rey \dT Q = \ &
  \frac{5}{2}\frac{\sin(\alpha{-}\theta)}{\sin\alpha}F{-}\frac{5}{2}\frac{Q}{F^2}
  {-}\frac{5}{2}\delta\frac{\cos(\alpha{-}\theta)}{\sin\alpha}\dX F F
  {-}\frac{15}{16}\delta\frac{\sin(\alpha{-}\theta)}{\sin\alpha}\dX\theta F^2 \notag\\
  &+\frac{5}{2}\ibond(\delta \dX[3]F{-}\zeta \dX K)F {-}\frac{17}{7} \delta\rey
  \frac{Q}{F} \dX Q +\frac{9}{7}
  \delta \rey \frac{Q^2}{F^2}\dX F+\frac{9}{2}\delta^2 \dX[2]Q \nonumber \\
  &+\frac{45}{16}\delta\zeta K \frac{Q}{F}+4\delta^2 \frac{Q}{F^2}(\dX
  F)^2{-}6\delta^2\frac{Q}{F}\dX[2]F {-}
  \frac{9}{2}\delta^2\frac{1}{F}\dX Q \dX F  \nonumber \\
  &+\delta^2\rey^2\biggl( {-}\frac{1}{210}\pa_{XT}\! Q Q F{-}\frac{1}{105} \dT Q
    \dX Q F{-}\frac{1}{21}(\dX Q)^2 Q
{-}\frac{1}{70}\dX[2]Q Q^2{+}
    \frac{1}{70}\frac{Q^2}{F}\dX Q \dX F \biggr)\notag
\end{align}
where we made use of \eqref{FT} to eliminate time derivatives of $F$. As there
are still time derivatives of $Q$ on the right-hand side this is not yet an
explicit evolution equation for $Q$.  However, from \eqref{3_28} we know that
$Q = F^3 + \ord(\eps)$, which leads to $ \dT Q = 3F^2 \dT F+\ord(\eps)
=-3\frac{Q}{F} \dX Q+\ord(\eps)$. Together with \reff{FT} this gives the
evolution system for $(F,Q)$, namely
\begin{align}
  \dT F =\ & -(1{-}\delta\zeta KF)\dX Q, \label{FT2}\\
  \delta\rey \dT Q = \ &
  \frac{5}{2}\frac{\sin(\alpha{-}\theta)}{\sin\alpha}F-\frac{5}{2}\frac{Q}{F^2}
  -\frac{5}{2}\delta\frac{\cos(\alpha{-}\theta)}{\sin\alpha}\dX F F
  -\frac{15}{16}\delta\frac{\sin(\alpha{-}\theta)}{\sin\alpha}\dX\theta F^2 \notag\\
  &+\frac{5}{2}\ibond(\delta \dX[3]F-\zeta \dX K)F -\frac{17}{7} \delta\rey
  \frac{Q}{F} \dX Q +\frac{9}{7}
  \delta \rey \frac{Q^2}{F^2}\dX F+\frac{9}{2}\delta^2 \dX[2]Q \nonumber \\
  &{+}\frac{45}{16}\delta\zeta K \frac{Q}{F}{+}4\delta^2 \frac{Q}{F^2}(\dX
  F)^2{-}6\delta^2\frac{Q}{F}\dX[2]F {-} \frac{9}{2}\delta^2\frac{1}{F}\dX Q \dX F
  {-}\frac{1}{210}\delta^2\rey^2(\dX Q)^2 Q. \label{QT}
\end{align}
If we set the waviness $\zeta = 0$ we obtain a system which is up to scaling
the same as the non-regularized WRIBL in \cite{Scheid_06}. That means that in case
of a flat bottom our one-step method is indeed equivalent to the Galerkin
method with universal polynomials and subsequent simplification. Thus for
$\zeta = 0$ our WRIBL is consistent with the Benney equation and
predicts the correct critical Reynolds number $\Rc$. In the next section we
will check the consistency for $\zeta > 0$ before we will regularize the equation in Section \ref{sec_pade}.

\section{Consistency}
\label{sec_cons}

The basic assumption throughout this paper is that $F$ is of order $1$ while
$U,W$ and $P$ can be expressed in powers of $\eps$ as stated in
\eqref{3_1}. In Section 3 this allowed us to solve the Navier--Stokes equations
asymptotically, which was used in 
Section 4 to derive the evolution equation \eqref{FT} for $F$ depending on the
flow rate $Q$. The natural approach to achieve a scalar equation is now to
plug into \eqref{FT} the expansion
\begin{equation}
  Q = Q_0+\eps Q_1+\eps^2 Q_2+\ord(\eps^3) = \int_0^F U_0 dZ+\eps \int_0^F U_1 dZ+\eps^2 \int_0^F U_2
  dZ+\ord(\eps^3).\notag
\end{equation}
We call the resulting equation Benney equation for wavy bottoms. In
\eqref{3_28} and \eqref{3_29} we have already calculated the zeroth and first
order components $Q_0$ and $Q_1$.  Consistency now means the following: In the
evolution equation \eqref{QT} for $Q$ we formally replace $Q$ by an enslaved
version $\QI$ with the expansion
\begin{equation}
  \QI = \QI_0+\eps \QI_1+\eps^2 \QI_2+\ord(\eps^3). \label{5_1}
\end{equation}
It is remarkable that $-\frac{5}{2}\frac{Q}{F^2}$ is the only $\ord(1)$-term
in \eqref{QT} which contains $Q$. Thus we obtain a set of linear algebraic
equations for $\QI_0, \QI_1, \QI_2$ which can be solved easily. By plugging
$\QI$ into \eqref{FT} we obtain a second scalar evolution equation for $F$. We
call our WRIBL consistent if this approach yields the Benney
equation for wavy bottoms.

To derive the Benney equation for wavy bottoms by a long
wave expansion of the Navier--Stokes equations and the associated boundary
conditions \eqref{NS1}--\eqref{NS6} we continue  as in 
\eqref{3_17}, \eqref{3_19}. 
 At $\ord(\eps^2)$ we obtain $U_2$. As this is
rather lengthy we refer to Appendix \ref{sec_u2} and state
here only the integrated version, namely
\begin{alignat}{2}
  \eps^2 Q_2 & = & \ & \eps^2 \int_0^FU_2 dZ \nonumber \\
  & = & & \frac{12}{7}\delta^2\rey^2\dX[2]F F^{10} + \frac{381}{35}\delta^2\rey^2(\dX F)^2 F^9 + \frac{10}{7}\delta^2 \rey(\ibond\dX[4]F - \cot\alpha \dX[2]F)F^7\nonumber\\
  & & & - \frac{8}{35}\delta\zeta\rey(\ibond \dX[2]K_0 + \cot\alpha\dX\theta_1)F^7 + \frac{12}{5}\delta\rey(3\delta\ibond(\dX[2]F)^2 - 2\delta\cot\alpha(\dX F)^2 \nonumber \\
  & & &{+} 5\delta\ibond\dX[3] F \dX F  {-} \zeta(\ibond\dX K_0 {+} \cot\alpha \ \theta_1) \dX F)F^6 {+} \frac{72}{5}\delta^2\rey\ibond(\dX F)^2\dX[2]F F^5 \nonumber \\
  & & &{+} \frac{9}{8}\delta\zeta K_0 F^4 {-}\frac{3}{8} \delta\zeta \dX \theta_1
  F^4 {+} 3\delta^2\dX[2]F F^4 {-} \delta\zeta\theta_1 \dX F F^3 {+} 7\delta^2(\dX
  F)^2 F^3 {-}\frac{1}{2} \zeta^2\theta_1^2 F^3.\label{q2}
\end{alignat}
Replacing $Q$ in \eqref{FT2} by $Q_0+\eps Q_1+ \eps^2 Q_2$ yields the Benney
equation for wavy bottoms.

Now we use \reff{QT} to derive a scalar model. Plugging \eqref{5_1} into
\eqref{QT} yields at $\ord(1)$:
\begin{equation}
  \frac{5}{2}F-\frac{5}{2}\frac{\QI_0}{F^2} = 0 \ \Leftrightarrow \ \QI_0 = F^3.
\end{equation}
At first order we get
\begin{align}
  \delta\rey \dT \QI_0 = & {-}\frac{5}{2}\zeta\cot\alpha \ \theta_1 F
  {-}\frac{5}{2}\eps\frac{\QI_1}{F^2}{-}\frac{5}{2}\delta\cot\alpha \dX F F
  {+}\frac{5}{2}\delta\ibond\dX[3]F F{-}\frac{5}{2}\zeta\ibond\dX K_0 F\nonumber\\
  &-\frac{17}{7}\delta\rey \frac{\QI_0}{F}\dX \QI_0 +
  \frac{9}{7}\delta\rey\frac{(\QI_0)^2}{F^2} \dX F.
\end{align}
By applying $\dT \QI_0 = 3 F^2 \dT F =-3F^2\dX \QI_0 {+} \ord(\eps) = -9\dX F
F^4 {+} \ord(\eps)$ this equation can be solved for $\QI_1$, which yields 
\begin{equation}
  \eps\QI_1 = \left(\frac{6}{5} \delta\rey \dX F F^3{-}\zeta(\ibond \dX K_0 
{+} \cot\alpha \ \theta_1){-}\delta\cot\alpha \ \dX F{+}\delta \ibond \dX[3] F\right) F^3.
\end{equation}
Comparing these results with \eqref{3_28}, \eqref{3_29} we already see that
$Q$ and $\QI$ match at zeroth and first order. In order to calculate $\QI_2$
we solve \eqref{QT} at $\ord(\eps^2)$. As this is somehow
elaborate and does not give any new insight we state here only the result,
i.e.~$\QI_2 = Q_2$ as expected.
As both the long wave expansion and the WRIBL approach yield the same
expansion of $Q$, the scalar evolution equations are in both cases
the same. Therefore our WRIBL is consistent with the Benney equation
also for $\zeta > 0$.

\section{Regularization}
\label{sec_pade}

With the WRIBL \eqref{FT2}, \eqref{QT} we now have a second-order
model for film flow over wavy bottoms which is consistent with the
according Benney equation and reproduces in the limit of a flat
incline the correct critical Reynolds number $\Rc$. In order to
achieve consistency the basic idea of the one-step Galerkin method was
to use as test and ansatz function a velocity profile which is a
solution of the expanded Navier--Stokes equations
\eqref{NS1}--\eqref{NS6} also in the time dependent case. Therefore in
\eqref{3_19} the first-order component $U_1$ in particular contains
the time derivative $\dT F$ which is substituted by the zeroth-order
identity $\dT F = - 3\dX F F^2$. In contrast to setting $\dT F = 0$ in
the velocity profile this procedure leads to the additional term
$-\frac{1}{210} \delta^2\rey^2(\dX Q)^2 Q$ in the WRIBL \eqref{FT2},
\eqref{QT} which turned out to be necessary for
consistency. 

However, over flat bottom it is known that a pure 
asymptotic expansion
approach with the above substitution of $\dT F$ can lead to an
unphysical behaviour if the Reynolds number exceeds a certain value
$\rey_0$ not far beyond $\Rc$.
In \cite{Pumir_83} one-hump solitary wave
solutions of a scalar Benney-like equation
for flat inclines are considered. 
According to the bifurcation diagram
\cite[Fig.~5]{Pumir_83} such homoclinic orbits are only found if the Reynolds
number is close to the instability threshold, i.e. $\Rc < \rey <
\rey_0$. However, in \cite{Salamon94}, where the two-dimensional
Navier--Stokes equations were solved by a finite-element method, such
a limit $\rey_0$ was not obtained. Thus the asymptotic expansion
equation used in \cite{Pumir_83} appears to be valid only 
if $\rey$ is not far beyond $\Rc$, and shows non-physical 
behaviour if $\rey$
exceeds a limiting value $\rey_0$. This deficiency appears to be
closely related to finite-time blow-up solutions in the scalar Benney
equation. 

For flat vertical walls it was shown in \cite{Scheid_06} using homoclinic 
continuation that such a
limitation also occurs for the second-order WRIBL, i.e., the
branch of homoclinic orbits again turns back if the Reynolds number
becomes too large, see \cite[Fig.~1]{Scheid_06}. However, if the
inertia correction term, which corresponds to $-\frac{1}{210}
\delta^2\rey^2(\dX Q)^2 Q$ in our notation, is neglected this
non-physical loss of solitary waves ceases. 
At least for small $\zeta$ and otherwise similar parameters 
as in \cite{Scheid_06} we must expect similar problems with our model.

In \cite{Scheid_06} a Pad\'{e}-like approximant technique is used to
regularize the WRIBL in case of a flat incline, see also 
\cite{Ooshida99} for the case of a scalar surface equation. The main idea is to
remove the dangerous second-order inertia terms by multiplying the
residual equation for $\dt Q$ with a suitable regularization factor
$S$. This procedure preserves the degree of consistency since the
second-order inertia terms are still implicitly included. This becomes
clear if one applies the zeroth-order identity $Q=F^3$ to $S$ which
yields the original non-regularized WRIBL. 
Homoclinic continuation now yields 
solitary wave solutions for the regularized model with no 
non-physical behaviour for $\rey > \Rc$ \cite{Scheid_06}. More
precisely, for a wide regime of unstable Reynolds numbers 
solitary wave solutions are found, with amplitudes only slightly 
smaller than those
obtained by numerics for the Navier--Stokes equations, 
in contrast to the regularization in \cite{Ooshida99}. 

For the undulated bottom we again closely follow \cite{Scheid_06}. 
First, we split \eqref{QT} into three parts, namely
\begin{equation}
  \res_1 := \delta\rey(-\dT Q -\frac{17}{7}\frac{Q}{F} \dX Q 
+\frac{9}{7}\frac{Q^2}{F^2}\dX F)
\quad\text{and}\quad \res_2 := -\frac{1}{210}(\delta\rey)^2(\dX Q)^2 Q
  \label{r2}
\end{equation}
containing the inertia terms with leading $\delta\rey$ and
$(\delta\rey)^2$, respectively, and the rest
\begin{align*}
  \res_0 := &\frac{5}{2}\frac{\sin(\alpha{-}\theta)}{\sin\alpha}F
{-}\frac{5}{2}\frac{Q}{F^2} {-}\frac{5}{2}\delta
\frac{\cos(\alpha{-}\theta)}{\sin\alpha}\dX F F{-}\frac{15}{16}\delta\frac{\sin(\alpha{-}\theta)}{\sin\alpha}\dX\theta F^2 
{+}\frac{5}{2}\ibond(\delta \dX[3]F{-}\zeta \dX K)F \\
  &+\frac{9}{2}\delta^2 \dX[2]Q {+}\frac{45}{16}\delta\zeta K
  \frac{Q}{F}{+}4\delta^2 \frac{Q}{F^2}(\dX
  F)^2{-}6\delta^2\frac{Q}{F}\dX[2]F {-}
  \frac{9}{2}\delta^2\frac{1}{F}\dX Q \dX F.
\end{align*}
The $\dT Q$-equation \eqref{QT} now reads $\res_0 + \res_1 +
\res_2 = 0$, and using again  $Q=F^3$ we see that $\res_2\sim (\dX F)^2 F^7$ 
is highly nonlinear. 
 The aim is to get rid of the potentially 
dangerous term $\res_2$
without loosing the degree of consistency. Therefore, {\em if we
enslave again $Q$ by $F$} as in Section \ref{sec_cons}, no term
up to $\ord(\eps^2)$ should be deleted or added. This is ensured,
e.g., if we multiply the residual equation by a regularization factor
$S$ which can depend on $F,Q$ and their derivatives. This yields $S
\res_0 + S(\res_1 + \res_2) =0$, and we are done if $S$ fulfills
\begin{equation}
  S(\res_1 + \res_2) = \res_1+\CO(\eps^3).
  \label{r3}
\end{equation}
This ansatz leads to the function
\begin{equation}
  S = \left( 1 + \frac{\res_2}{\res_1} \right)^{-1}.
\end{equation}
Plugging the zeroth-order identity $Q=F^3$ into \eqref{r2} yields
\begin{equation*}
  \res_1 = 3\delta\rey\dX F F^4 + \ord(\eps^2), \quad 
\res_2 = -\frac{3}{70}(\delta\rey)^2(\dX F)^2 F^7 + \ord(\eps^3), 
\end{equation*}
and thus, using again  $Q=F^3$, 
\begin{equation}
  \tilde S := \left(1 - \frac{1}{70}\delta\rey Q \dX F\right)^{-1} 
= S + \ord(\eps^2). 
\end{equation}
Then \eqref{r3} leads to $\tilde S(\res_1 + \res_2) = \res_1 +
\ord(\eps^3)$, and multiplying $\res_0 + \res_1 +
\res_2 = 0$ by $\tilde S$ finally yields the ``regularized'' equation
$\tilde S \res_0 + \res_1 = \ord(\eps^3)$. In summary, the regularized 
version (rWRIBL) of
the weighted residual integral boundary layer equation reads
\begin{align}
  \dT F =\ & -(1{-}\delta\zeta KF)\dX Q, \label{FTR}\\
  \delta\rey \dT Q = \ &
  -\frac{17}{7} \delta\rey \frac{Q}{F} \dX Q +\frac{9}{7} \delta \rey \frac{Q^2}{F^2}\dX F +\left(\frac{5}{2}\frac{\sin(\alpha{-}\theta)}{\sin\alpha}F-\frac{5}{2}\frac{Q}{F^2}-\frac{5}{2}\delta\frac{\cos(\alpha{-}\theta)}{\sin\alpha}\dX F F\right. \nonumber \\
  & -\frac{15}{16}\delta\frac{\sin(\alpha{-}\theta)}{\sin\alpha}\dX\theta F^2 +\frac{5}{2}\ibond(\delta \dX[3]F-\zeta \dX K)F +\frac{9}{2}\delta^2 \dX[2]Q {+}\frac{45}{16}\delta\zeta K \frac{Q}{F}\nonumber \\
  & \left.{+}4\delta^2 \frac{Q}{F^2}(\dX
    F)^2{-}6\delta^2\frac{Q}{F}\dX[2]F {-}
    \frac{9}{2}\delta^2\frac{1}{F}\dX Q \dX F\right)\left(1 -
    \frac{1}{70}\delta\rey Q \dX F\right)^{-1}. \label{QTR}
\end{align}

It is not easy to assess the value of this regularization. First, for 
flat bottom we numerically confirmed 
the loss of the one-hump solitary waves for the WRIBL \reff{FT2}, \reff{QT} 
in a certain interval $[\rey_0,\rey_1]$ of $\rey>\Rc$ and its regain 
for the rWRIBL \reff{FTR}, \reff{QTR}. 
However, here we use direct numerical simulations (see Section \ref{sec_num} 
for details), instead of homoclinic 
continuation in \cite{Scheid_06}, which is not possible for $\zeta>0$, or 
in any case is much more involved since the solitary waves then do not 
decay to a constant state but to spatially periodic solutions. 
In these direct numerical simulations we find 
that the interval $[\rey_0,\rey_1]$ is 
typically rather narrow, shrinks quickly with increasing $\zeta>0$ 
and vanishes for $\zeta$ 
greater some $\zeta_0$ which depends on the other parameters. 
Also, the loss of solitary waves in $[\rey_0,\rey_1]$ is not related 
to blow-up of solutions: instead, small amplitude irregular patterns 
appear in this interval. This might indicate a transition between two 
different branches of solitary waves for $\rey<\rey_0$ and $\rey>\rey_1$, 
or some other more complicated structure in the background. 

\begin{figure}[!ht]
\begin{center}
 \begin{tabular}[h]{m{0.1cm}@{}>{\centering\arraybackslash}m{6.1cm}}
    \begin{sideways}
      \footnotesize Dimensionless thickness $F$ 
    \end{sideways}
    & \ig[width=60mm]{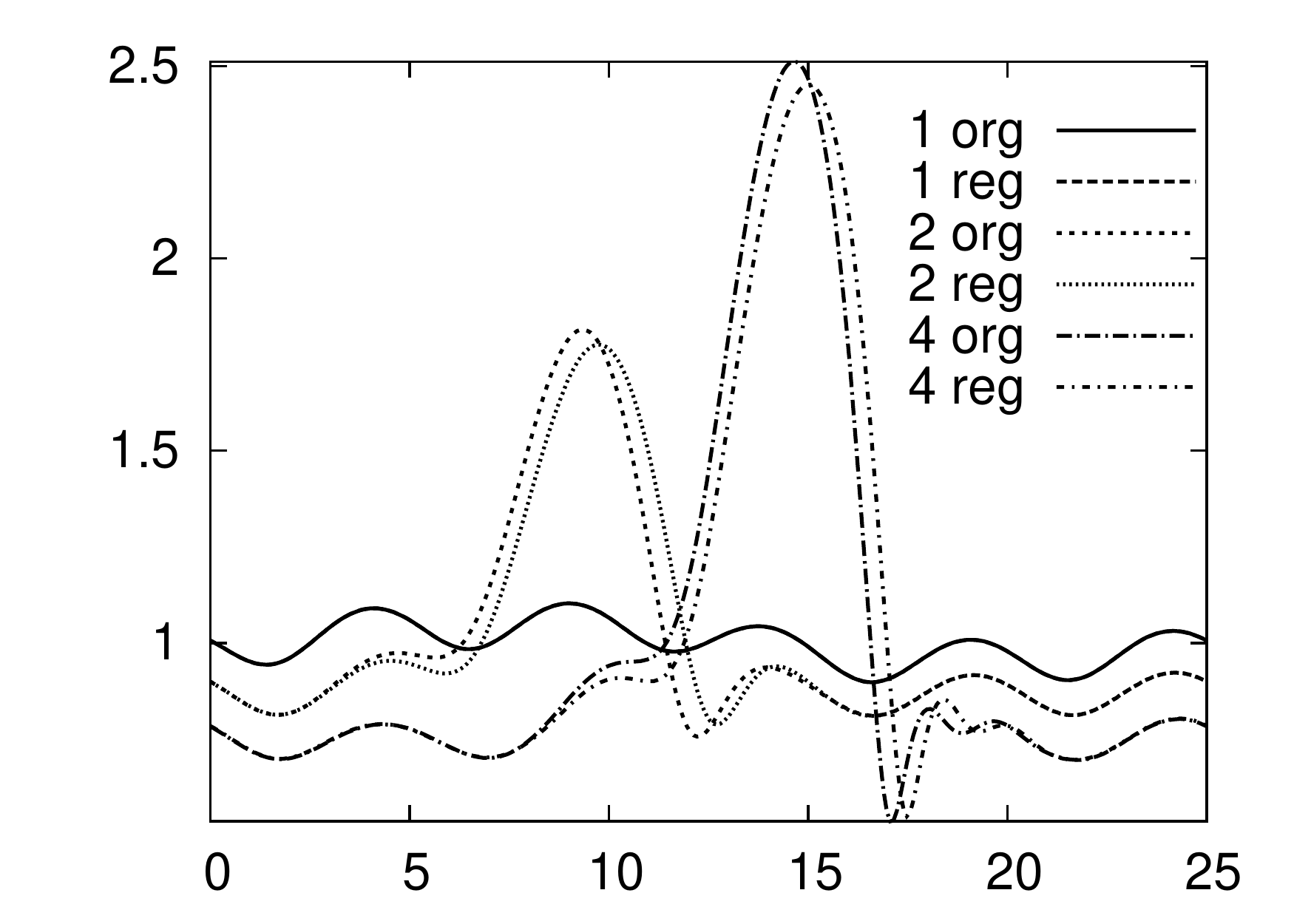} \vspace{-0.2cm} \\
    \small (a) &  \footnotesize $\xhat$ [mm]
 \end{tabular}
 \begin{tabular}[h]{m{0.1cm}@{}>{\centering\arraybackslash}m{6.1cm}}
    \begin{sideways}
      \footnotesize Max. amplitude of $F$ 
    \end{sideways}
    & \ig[width=60mm]{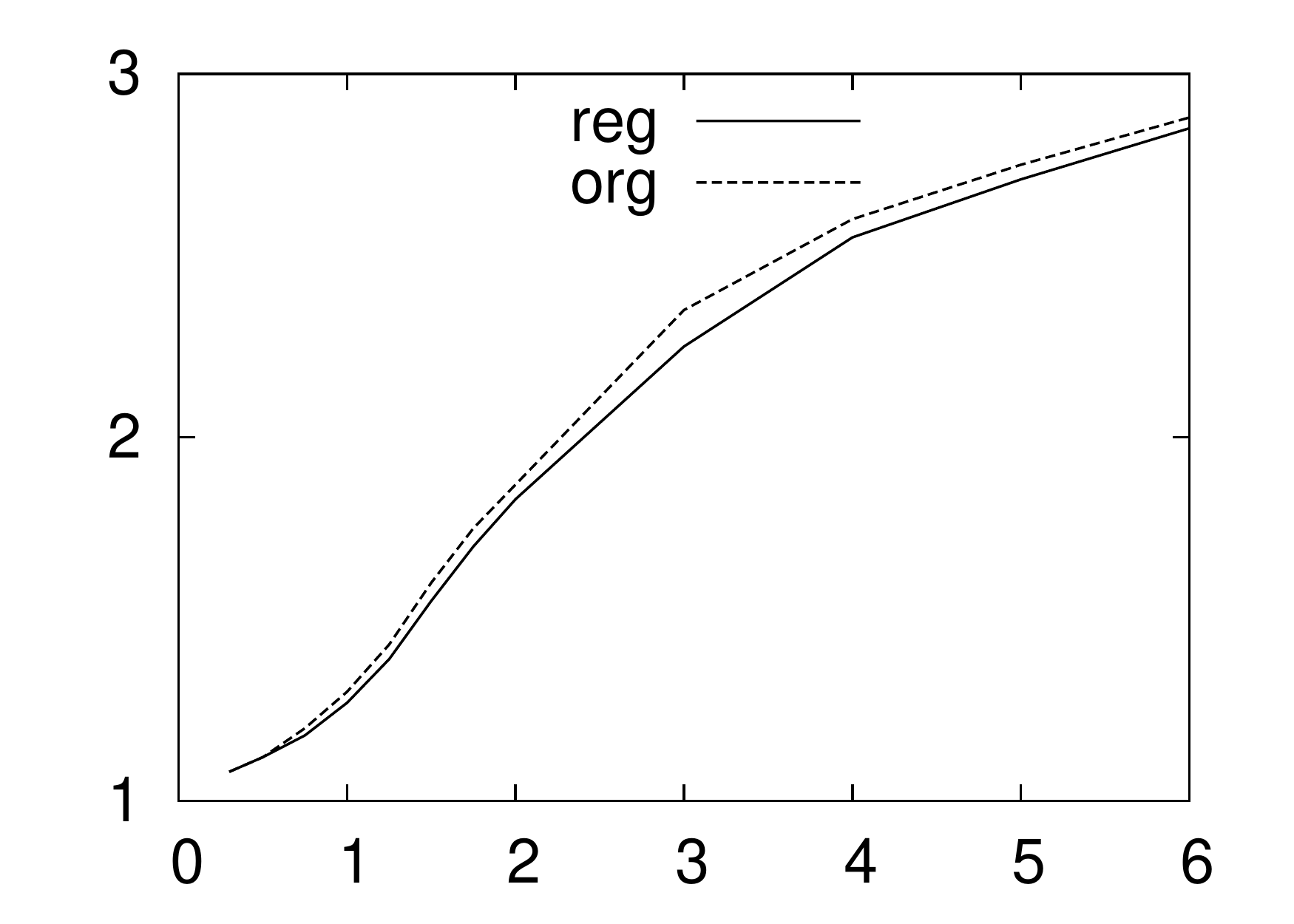} \vspace{-0.2cm} \\
    \small (b) &  \footnotesize $\rey$
	\end{tabular}
  \caption{Comparison of the WRIBL with the regularized version 
rWRIBL. $\al=90^\circ, \del = 0.3, \zeta = 0.05, \ibond=3.32$ (comparable to 
\cite[Fig.~1]{Scheid_06}); $\lambdahat=5 \, \text{mm}$, 5 bottom waves, 
$\rey$ as indicated, and {\tt org} and {\tt reg} stand 
for the original WRIBL and the regularized version rWRIBL. 
(a) shows snapshots of the dimensionless film thickness 
$F(\xhat)$ with $\spr{F}=1$, and (b) the maximal amplitude of $F$ 
extracted from one time period of well 
converged one-hump solitary waves. 
 For these parameters, solitary waves 
of both the WRIBL and the rWRIBL are found for all $\rey\in (0.3,\rey_2)$ 
with $\rey_2\approx 12$, where (for the used discretization $n=400$) 
the numerics fail for the rWRIBL due to $F\to 0$ pointwise. 
Generically, the solitary waves for 
the rWRIBL have slightly smaller amplitude. }
  \label{rfig}
\end{center}
\end{figure}

To illustrate the effect of the regularization, Fig.~\ref{rfig} 
shows (in advance of \S\ref{sec_num}) some differences between the 
rWRIBL and the WRIBL for a parameter set 
for which there is no interval 
 $[\rey_0,\rey_1]$ where the WRIBL does not have solitary wave solutions 
in direct numerical simulations. 
In general, 
these differences appear to be rather small, with the notable exception of 
the calculation of the critical Reynolds number $\Rc$ in 
Fig.~\ref{rc-fig} below, where the results for the rWRIBL are closer 
to available data. 

In general, in our simulations both the WRIBL and the rWRIBL did not show 
blow-up of solutions in parameter regimes of interest, 
but there appears to be one disadvantage of the rWRIBL: 
for some parameters, as $\rey$ becomes large the numerics for the 
rWRIBL fail more rapidly than those for the WRIBL. In particular, 
for the parameters in Fig.~\ref{rfig} we can follow one-hump 
solitary waves for the WRIBL up to $\rey\approx 90$ where these 
split up into two humps, while for the rWRIBL we obtain numerical 
failures due to $F\to 0$ pointwise for $\rey$ not far beyond 12. 
However, this is strongly related to the method of simulation, 
i.e., to the fact that 
$\spr{F}=1$ is imposed, and should not be considered as blow-up 
of solutions of the rWRIBL: for instance we can follow 
 one-hump solitary waves for the rWRIBL up to $\rey=21$ if 
we double the domain length in Fig.~\ref{rfig}. 
In summary, since we are more interested in the regime $\rey$ not 
too far from $\Rc$, where the rWRIBL gives results closer 
to available data than the WRIBL, below we focus 
on the rWRIBL for our numerical simulations.

\section{Numerical simulations}
\label{sec_num}
Though the rWRIBL \reff{FTR}, \reff{QTR} is much simpler than the Navier--Stokes
system \reff{eq1}--\reff{eq7}, it is still a quasilinear parabolic system,
with periodic coefficients.  Therefore, a first step to explore some of its 
stationary and non-stationary solutions are numerical simulations. 
For this we have set up a finite difference
method with periodic boundary conditions in space for 
both, the rWRIBL and the WRIBL.  To calculate 
stationary solutions  $(F,Q)_s$ we use a Newton method starting 
at constant $(F,Q)$ which corresponds to a Nusselt flow, 
which in contrast to the flat bottom case is {\em not} a stationary solution 
over wavy bottom. For the time dependent problem we may also use 
 constant $(F,Q)$ or perturbations of 
some $(F,Q)_s$  as initial data. We then use an implicit and 
adaptive time stepping.  Depending on the flow characteristics, 
the spatial discretization was on the order of 50 (Fig.~\ref{rc-fig})
to 400 (Fig.~\ref{ngf}) points per 
bottom wave. Numerical convergence was checked 
by refining the discretization without perceivable differences in the 
solutions. 

\subsection{Comparison with available data}\label{onum}
First we want to compare our results with available experimental and 
numerical data.  Therefore we 
have to somewhat relax the assumption used in the derivation 
of the WRIBL that $\rey$ and $\ibond$ are of order $1$
compared to $\zeta,\delta$ which are assumed to be small. 
However, similar relaxations often appear in the {\em application} 
of asymptotic expansions. In other words, one goal of the present 
section is to study how far the asymptotic expansion can take us. 
As said above, we focus on 
the rWRIBL since it gives slightly better comparison with available 
data. 

\begin{figure}
\begin{center}
  \begin{tabular}[h]{m{0.2cm}m{6cm}}
    \begin{sideways}
      \footnotesize Local film thickness [mm]
    \end{sideways}
    & \ig{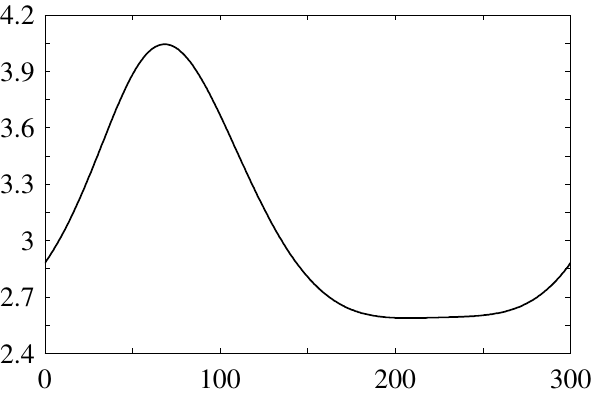} \vspace{-0.7cm} \\
    \small (a) & \centering \footnotesize $\xhat$ [mm]
  \end{tabular}
  \begin{tabular}[h]{m{0.2cm}m{6cm}}
    \begin{sideways}
      \footnotesize Local film thickness [mm]
    \end{sideways}
    & \ig{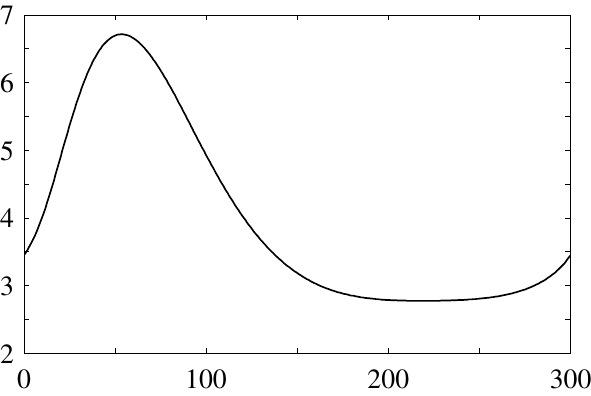} \vspace{-0.7cm} \\
    \small (b) & \centering \footnotesize $\xhat$ [mm]
  \end{tabular}
  \caption{Local film thickness for two different inclination angles. For
    comparison with \mbox{\cite[Fig.~3]{Acta_05}} it is measured not perpendicular to
    the bottom but to the main flow direction $\exhat$, see
    \mbox{Fig.~\ref{fig1}}. Parameters: $\rey = 0.0285, \zeta = 0.31$ and (a)
    $\alpha = 28^\circ, \delta = 0.059, \ibond = 2\times 10^{-3}$, \mbox{(b) $\alpha
    = 18.05^\circ,$} $\delta = 0.068, \ibond = 3\times 10^{-3}$.}
  \label{fig_surf}
\end{center}
\end{figure}

\begin{figure}
\begin{center}
  \begin{tabular}[h]{m{0.2cm}m{6cm}}
    \begin{sideways}
      \footnotesize Film thickness [mm]
    \end{sideways}
    & \ig[width=6cm]{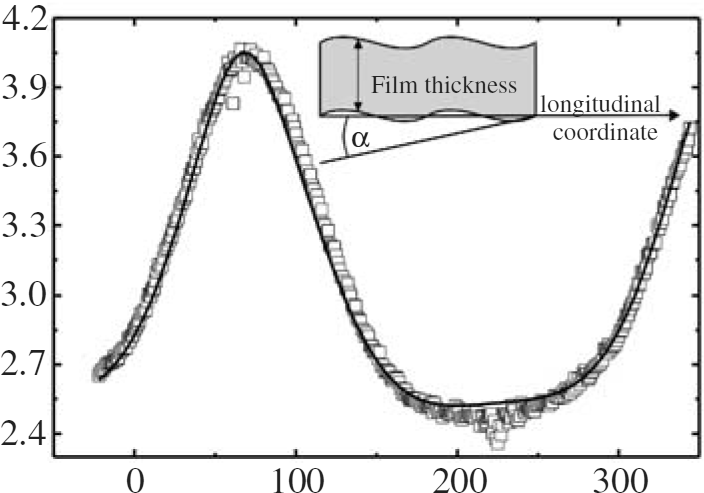} \vspace{-0.5cm} \\
    \small (a) & \centering \footnotesize Longitudinal coordinate [mm]
  \end{tabular}
  \begin{tabular}[h]{m{0.2cm}m{6cm}}
    \begin{sideways}
      \footnotesize Film thickness [mm]
    \end{sideways}
    & \ig[width=6cm]{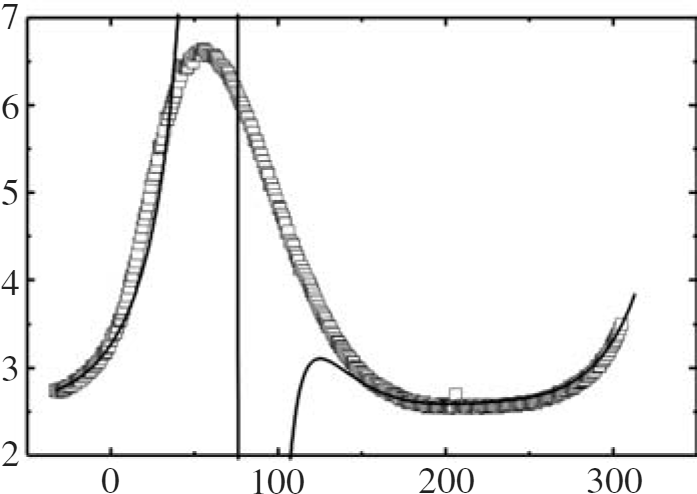} \vspace{-0.5cm} \\
    \small (b) & \centering \footnotesize Longitudinal coordinate [mm]
  \end{tabular}
  \caption{Experimental data for the parameters used in Fig.~\ref{fig_surf}. Reprint of \cite[Fig.~3]{Acta_05}, with permission from Springer Science+Business Media.}
  \label{fig_surf_orig}
\end{center}
\end{figure}

We first  simulate the stationary problem 
for fluid and geometry parameters taken from \cite{Acta_05}, namely 
$\nu = \unit{1110}{\milli\squaren\meter\per\second}, 
\rho = \unit{0.969}{\gram\per\centi\cubic\meter}$, 
$\sigma = \unit{20.4}{\milli\newton\per\meter}$. The bottom is a sine with 
wavelength $\lambdahat = 300 \, \text{mm}$, amplitude $\ahat
= 15 \, \text{mm}$ and trough and crest at $\xhat=0,\xhat=150$, respectively. 
Fig.~\ref{fig_surf} shows the resulting local film
thickness which is the distance of the free surface to the bottom contour
measured in $\ezhat$-direction, see Fig.~\ref{fig1}. As inclination 
angles we take (a) $\alpha =
28^\circ$, (b) $\alpha = 18.05^\circ$. Choosing the Reynolds number such that
the maximum local film thickness is the same as in \cite[Fig.~3]{Acta_05} we
obtain stationary solutions $(F,Q)_s$ 
which for the film height are in perfect agreement with experimental 
data, see Fig.~\ref{fig_surf_orig}. 

\begin{figure}[!ht]
\begin{center}
  \begin{tabular}[h]{@{}m{0.3cm}@{}m{11.5cm}@{}}
    \begin{sideways}\footnotesize $\hat{z}$ [mm]\end{sideways}
    & \ig{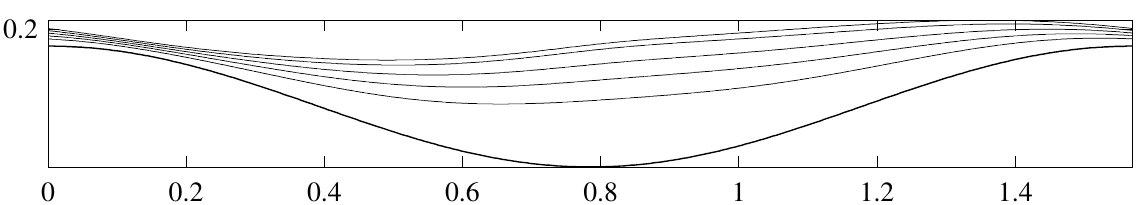} \\
      \vspace{-0.4cm}\small (a) & \vspace{-0.4cm}\centering \footnotesize $\xhat$ [mm]
 \end{tabular}\\ \vspace{0.2cm}
  \begin{tabular}[h]{@{}m{0.3cm}@{}m{11.5cm}@{}}
    \begin{sideways} \footnotesize $\hat{z}$ [mm] \end{sideways}
    & \ig{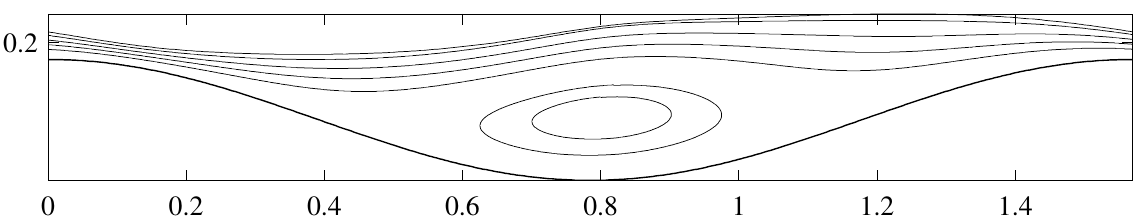} \\
      \vspace{-0.4cm}\small (b) & \vspace{-0.4cm}\centering \footnotesize $\xhat$ [mm]
  \end{tabular}
  \caption{Free surface and reconstructed flow field for stationary
    solutions of \eqref{FTR}, \eqref{QTR} for (a) $\rey=5$ and (b) $\rey=20$. The
    other parameters are $\delta=0.15$ respectively $\delta=0.24$, $\alpha = 90^\circ,
    \ibond=17.92, \zeta=0.35, \lambdahat=1.57 \, \text{mm}.$}
  \label{fig_prof2}
\end{center}
\end{figure}

\begin{figure}[!ht]
\begin{center}
  \begin{tabular}[h]{@{}m{0.4cm}@{}m{11.5cm}@{}}
    \begin{sideways}\footnotesize $\hat{z}$ [mm]\end{sideways}
    & \ig[width=116mm]{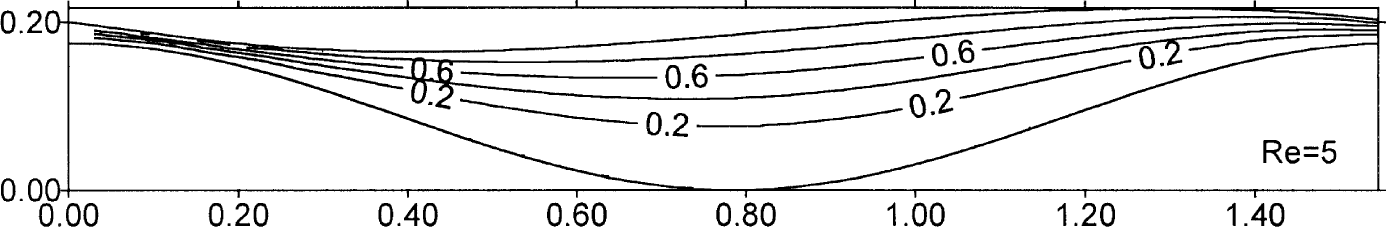} \\
      \vspace{-0.2cm}\small (a) & \vspace{-0.2cm}\centering \footnotesize Distance along gravity [mm]
 \end{tabular}\\ \vspace{0.2cm}
  \begin{tabular}[h]{@{}m{0.4cm}@{}m{11.5cm}@{}}
    \begin{sideways} \footnotesize $\hat{z}$ [mm] \end{sideways}
    & \ig[width=116mm]{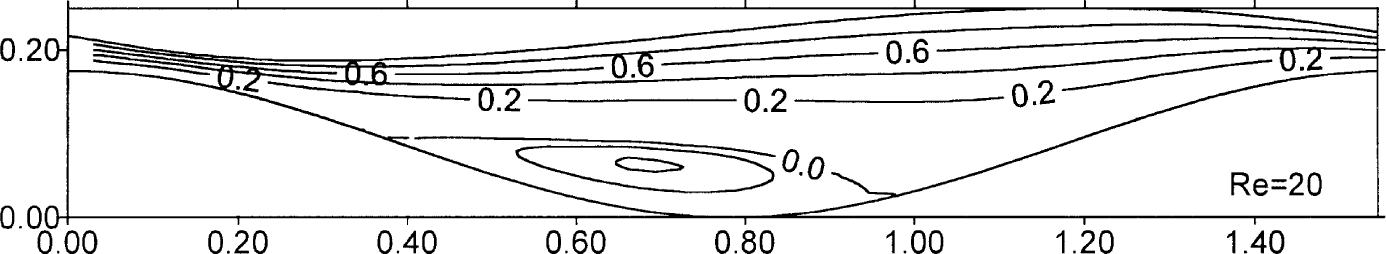} \\
      \vspace{-0.2cm}\small (b) & \vspace{-0.2cm}\centering \footnotesize Distance along gravity [mm]
  \end{tabular}
  \caption{Full Navier--Stokes numerics for the parameters used in Fig.~\ref{fig_prof2}. Reprint of \cite[Fig.~10]{Trifonov_98}, with permission from Elsevier.}
  \label{fig_prof2_orig}
\end{center}
\end{figure}

In order to explore wider regimes of parameters and to 
get more detailed comparison also with full Navier--Stokes numerics
we reconstruct the flow field using the second-order profile \eqref{Uprof2},
\eqref{Wprof2} derived in Appendix \ref{sec_u2} for sinusoidal bottoms.
Following \cite{Trifonov_98}, see also \cite{trif07a,trif07}, we 
simulate the flow of liquid nitrogen over a vertical sinusoidal
bottom with wavelength $\lambdahat = 1.57 \, \text{mm}$ and amplitude $\ahat =
0.0875 \, \text{mm}$. The fluid parameters are $\nu =
\unit{0.182}{\milli\squaren\meter\per\second}, \rho =
\unit{0.808}{\gram\per\centi\cubic\meter}$ 
and \mbox{$\sigma = \unit{8.87}{\milli\newton\per\meter}$} 
which yield an inverse Bond number $\ibond = 17.92$. 
As Reynolds numbers we choose $\rey = 5$ and $\rey =20$. 
Again we achieve free surface profiles which are in good agreement with 
the Navier--Stokes numerics in
\cite[Fig.~10]{Trifonov_98}, and also the flow fields are qualitatively 
and semi-quantitatively reproduced
correctly, see Fig.~\ref{fig_prof2} and \ref{fig_prof2_orig}. Namely, 
there occurs a 
recirculation zone of correct size in the 
trough of the bottom contour if the Reynolds number is increased. 

Above we calculated stationary solutions  $(F,Q)_s$ 
which, by analogy with the flat bottom case, 
must be expected to be unstable in the considered 
regime $(\alpha = 90^\circ)$, see 
also \cite{wa03,avb06,trif07a,trif07}. 
In the following we report on some numerical experiments 
to investigate the stability of stationary solutions and on some time 
dependent solutions in the unstable case. 
The standard approach to study the stability of $(F,Q)_s$ would be to 
calculate the spectrum of 
the linearization of \reff{FTR}, \reff{QTR} around $(F,Q)_s$, either 
numerically or 
analytically by expansion of first the stationary solution and 
then the eigenvalue problem in suitable small parameters.  
Eigenvalues of the linearization can then be calculated using 
Floquet theory. See \cite{trif07a} for a detailed parametric 
study of stability using this approach for an IBL, and 
\cite{trif07} for the full Navier--Stokes problem. 

Here, since we are mainly interested in the shape of non-stationary 
bifurcated solutions in case of instability, to determine stability 
of $(F,Q)_s$ we rather use a less systematic  {\it ad hoc} approach. 
We  numerically calculate $(F,Q)_s$ for various $\rey$, with
fluid and geometry parameters fixed.
Then,  on a domain with eight bottom undulations, we apply a 
localized perturbation, let the system run, and determine stability 
by growth or decay of the perturbations. This yields 
a critical Reynolds number $\Rc$ in terms of the remaining parameters. 

\begin{figure}[!ht]
  \centering
  \begin{tabular}[h]{m{0.3cm}m{8.4cm}}
    \begin{sideways}
      \footnotesize $\Rc$
    \end{sideways}
    & \hspace{-0.2cm}\ig{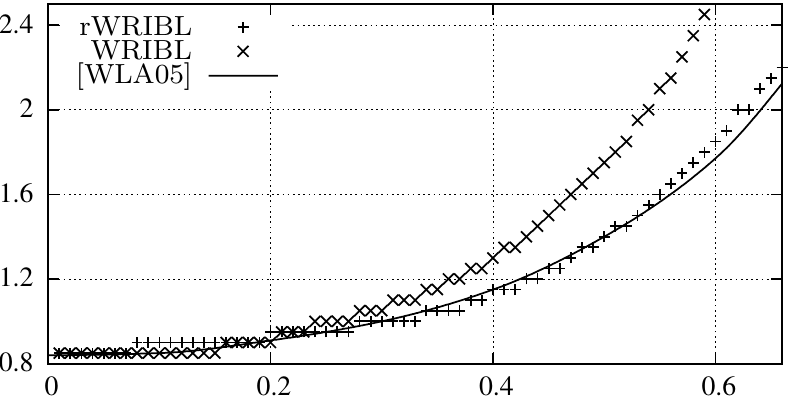} \vspace{-0.1cm}   \\
    \small & \centering \footnotesize $\zeta$
  \end{tabular}
  \caption{Critical Reynolds number $\Rc$ as a function of the waviness
    $\zeta$ for parameter set A from Table \ref{tab2}. Along the critical values $\del$
    varies from $\del = 0.035$ $(\rey=5/6)$ to $\del=0.048$ $(\rey=2.2)$. [WLA05] denotes $\Rc$ from 
\cite{Acta_05}, multiplied by 2/3 due to a different scaling. 
The critical Reynolds numbers were calculated with a tolerance of $\pm 0.05$. 
}
  \label{rc-fig}
\end{figure}

\begin{table}[!ht]
\begin{center}
  \begin{tabular}{|l|c|c|c|}
    \hline
     & A & B & C \\ \hline \hline
    $\rho [\unit{}{\gram\per\centi\cubic\meter}]$       & $0.969$ & $0.969$ & $1.00$ \\ \hline
    $\nu [\unit{}{\milli\squaren\meter\per\second}]$    & $24.1$  & $24.1$  & $1.00$ \\ \hline
    $\sigma [\unit{}{\milli\newton\per\meter}]$         & $20.0$  & $20.0$  & $70.0$ \\ \hline
    $\lambdahat [\unit{}{\milli\meter}]$                & $108$   & $108$   & $10.0$ \\ \hline
    $\alpha [^\circ]$                                   & $45$    & $10$    & $10$   \\ \hline\hline
    $\ibond$                                            & $0.01$  & $0.04$  & $16.2$ \\ \hline
  \end{tabular}
\caption{Parameters used to study stability of stationary solutions, with resulting inverse Bond numbers.}\label{tab2}
\end{center}
\end{table}
Again we first focus on non-dimensional parameters from \cite{Acta_05}, namely
$\alpha = 45^\circ$ and \mbox{$\ibond = 0.01$}, using the dimensional
 parameter set A from Table \ref{tab2}, and calculate $\Rc$ as function of 
$ \zeta$, see Fig.~\ref{rc-fig}. 
In agreement with \cite[Fig.~7]{Acta_05}, see also \cite{avb06,trif07}, we 
find that the wavy bottom strongly increases $\Rc$ compared to
the critical Reynolds number $5/6\cot\alpha$ over flat bottom. In particular, 
also the quantitative agreement with \cite[Fig.~7]{Acta_05} is very good. 
Here the most notable difference between the WRIBL and the rWRIBl occurs:  
$\Rc$ is somewhat larger for the WRIBL and hence the rWRIBL appears to be 
more accurate. 

\begin{figure}[!ht]
\begin{center}
  \begin{tabular}[h]{b{0.1cm}@{}b{5cm}@{}}
    \small (a) & \ig[width=50mm]{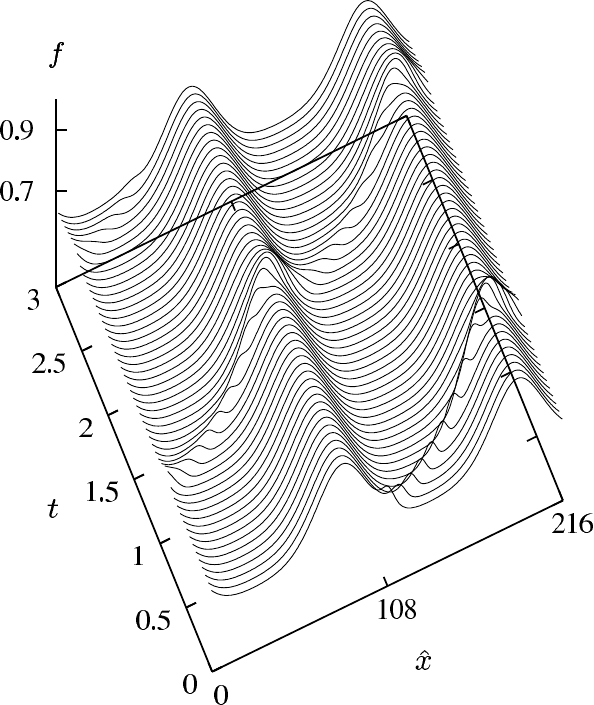}  
  \end{tabular}\hspace{0.6cm}
  \begin{tabular}[h]{b{0.1cm}@{}b{5cm}@{}}
    \small (b) & \ig[width=50mm]{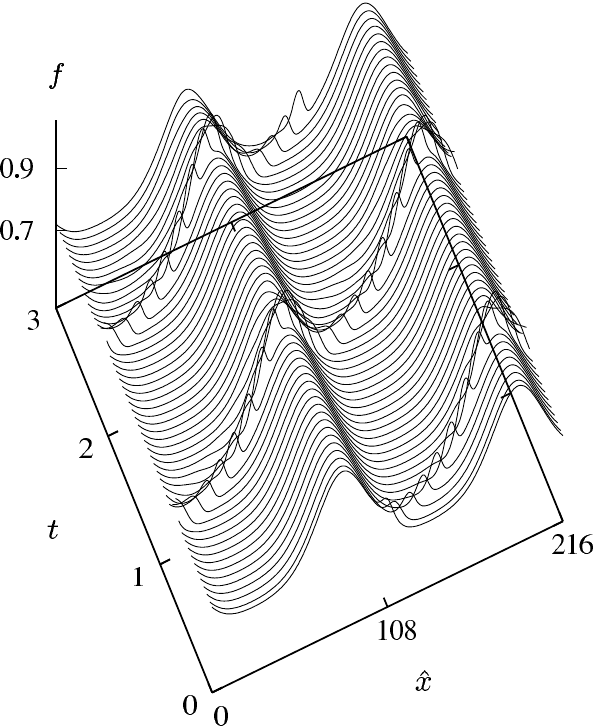}  
  \end{tabular}
  \caption{Numerical simulations in the sub- resp.~supercritical case 
for parameter 
set A from Table \ref{tab2} and $\zeta = 0.5$ which gives  $\Rc\approx 1.4$, cf.~Fig.~\ref{rc-fig}; two bottom waves with 
    periodic boundary conditions. (a),  $\rey=1.1$, $f(\xhat)$; 
for larger $t$ the
    solution relaxes to a stationary solution. (b) $\rey=1.6$, $f(\xhat)$;
    the solution is unstable and a traveling pulse evolves.}
  \label{dyn-fig}
\end{center}
\end{figure}

Figure \ref{dyn-fig} shows time dependent solutions, 
with $\zeta=0.5$ from Fig.~\ref{rc-fig}, but for graphical
reasons with only two bottom undulations. Over flat bottoms, 
for $\rey>\Rc$ the most prominent solutions are the 
(experimentally, numerically and
analytically well known) traveling pulse trains \cite{cd02}. Also over
wavy bottoms pulse like surface waves develop, and the effect of the
bottom waviness is a periodic modulation of the amplitude and speed of
the pulses: (a) shows the decay of a localized perturbation in the 
stable case, while (b) shows the emergence of a pulse in the unstable case. 

\subsection{Some new predictions}\label{nnum}
\begin{figure}[!htbp]
  \centering
  \begin{tabular}[h]{m{0.5cm}@{}m{8cm}}
    \begin{sideways}
      \footnotesize $\Rc$
    \end{sideways}
    & \ig{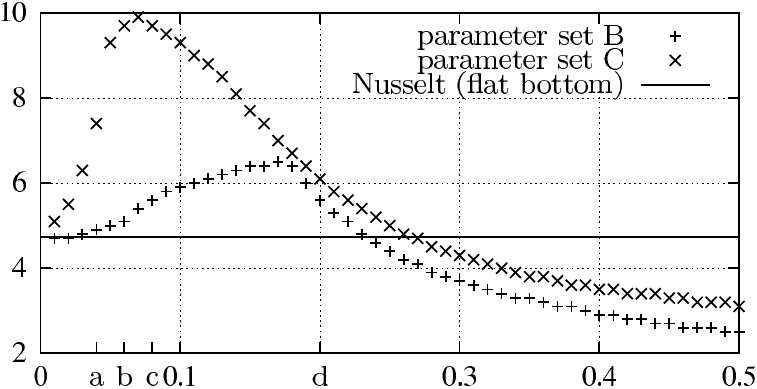} \vspace{-0.1cm}   \\
    \small & \centering \footnotesize $\zeta$
  \end{tabular}
  \caption{Critical Reynolds number $\Rc$ as a function of the
    waviness $\zeta$ for parameter sets B ($\ibond = 0.04$) and C
    ($\ibond = 16.2$) from Table \ref{tab2}; for parameter set C the
    letters a--d indicate the values of $\zeta$ used for the time
    dependent plots in Fig.~\ref{rc-fig-45}.}
  \label{rc-fig-23}
\end{figure}

The numerics in \S\ref{onum} have shown that \reff{FTR}, \reff{QTR} reproduces
known phenomena qualitatively and quantitatively, in particular the
appearance of eddies in troughs of the bottom for larger $\zeta$, and
the occurrence of a long wave instability when the Reynolds number
exceeds a critical value $\Rc$ as well as the increase of $\Rc$ with
$\zeta$. Next we consider a lower inclination angle for which we again
investigate the stability of stationary solutions by the method
specified above. Taking the same fluid parameters as in parameter set
A but with $\alpha = 10^\circ$ we get the critical values in
Fig. \ref{rc-fig-23} denoted by parameter set B. In contrast to
Fig.~\ref{rc-fig} the critical Reynolds numbers are no longer
increasing monotonously but reach a maximum at $\zeta \approx
0.17$. For larger values of the bottom waviness $\Rc$ decreases, and
for $\zeta > 0.23$ it becomes less than the critical Reynolds number
$5/6 \cot \alpha$ for flat bottom.

\begin{figure}[!ht]
  \begin{center} \small
    \begin{tabular}[h]{b{5.8cm}b{5.8cm}}
      \ig[width=55mm]{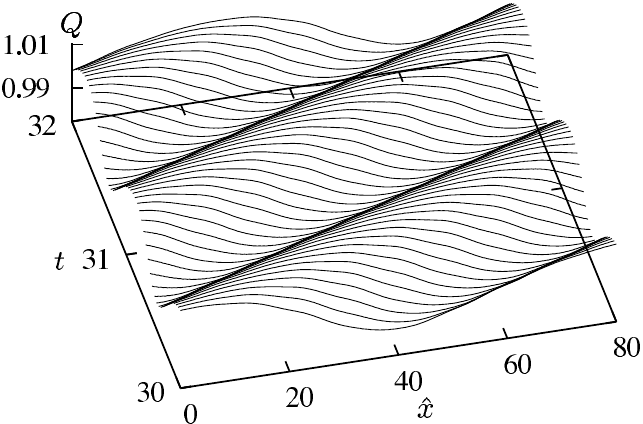} \vspace{-5mm} &  
\ig[width=55mm]{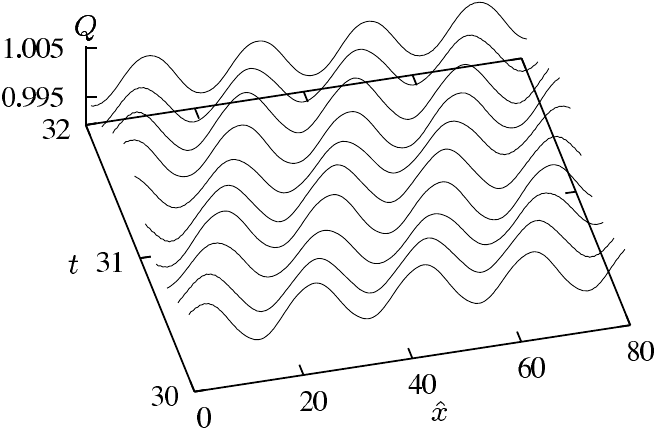} \vspace{-5mm} \\
      (a) & (b) \\
      \ig[width=55mm]{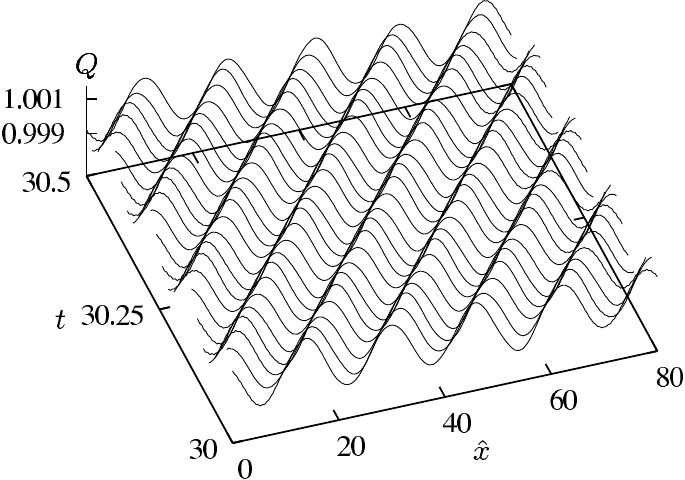} \vspace{-5mm} &  
\ig[width=55mm]{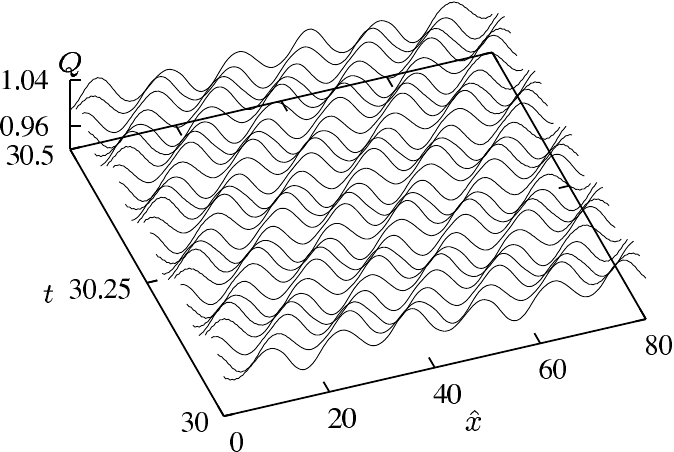} \vspace{-5mm} \\
      (c) & (d)
    \end{tabular}
    \caption{Time dependent simulations for parameter set C, eight
      bottom waves, plots of the flow rate $Q$. (a) $\zeta=0.04,
      \rey=7.4$ (long wave instability), (b) $\zeta=0.06, \rey=9.7$
      (short wave instability, four waves), (c) $\zeta=0.08, \rey=9.7$
      (short wave instability, five waves), (d) $\zeta=0.2, \rey=6.1$
      (short wave instability, six waves).}
    \label{rc-fig-45}
  \end{center}
\end{figure}

Next we increase the inverse Bond number by choosing $\lambdahat =
10\unit{}{\milli\meter}$ and the fluid parameters of water, see
parameter set C in Table \ref{tab2} and the resulting critical values
in Fig.~\ref{rc-fig-23}. The dependence on $\zeta$ turns out to be
more pronounced than in Fig.~\ref{rc-fig}. 
Figure \ref{rc-fig-45} shows related
time dependent solutions for some supercritical values. For small $\zeta$,
e.g.~$\zeta=0.04$ in (a), the instability is long wave (pulses), but
for $\zeta = 0.06$ in (b) the perturbation evolves into a finite
wavelength pattern. Thus, for $\alpha =10^\circ$ and $\zeta$ larger
than a critical value $\zeta_0 \lesssim 0.06$ there appears a 
finite wave number instability, and the wave number increases as the bottom
waviness becomes larger, see Fig.~\ref{rc-fig-45} (c)--(d). 
Since also the amplitudes of these patterns are very small we 
conclude on a phenomenological basis that Fig.~\ref{rc-fig-45} (b)--(d) 
shows short wave instabilities, where, however, 
the following remarks apply. 

The linearization of \reff{FTR}, \reff{QTR} around some $(Q,F)_s$ 
always has a Floquet exponent $\mu_1(0)=0$ from conservation of mass. 
In other words,  $\mu_1(0)=0$ since we have a family of stationary 
 solutions $(Q,F)_s$ parameterized by the total mass 
 $M=\int_0^{2\pi} F(1+\frac{1}{2}\delta\zeta KF)d\Xhat$. 
If $K\mapsto \mu(K)$ is a parameterization 
of the Floquet exponents of the linearization by wave number, then short 
wave instability in a strict sense means that unstable Floquet modes 
appear only in an interval $\pm K\in(K_1,K_2)$ with $K_1>0$. 
However, a finite wave number instability  
may also be due to a side band (i.e.~long wave) instability, 
that is, a branch $\mu(K)$ of unstable 
Floquet exponents with Re$\mu(K)=c_2K^2-c_4K^4+\CO(K^6)$ 
with $c_2,c_4>0$, which up to order $K^4$ gives $K_c=\sqrt{\frac{c_2}{2c_4}}$ 
as the most unstable wave number. To distinguish this from a short 
wave instability one should actually calculate the spectrum. 
However, we take (b)--(d) as strong hints for a short wave instability, 
since if (b)--(d) were due to side band instabilities we would expect 
larger amplitudes.  In any case, to distinguish (a) from (b)--(d) 
we may call the latter short wave instabilities in a 
phenomenological sense.

Finally, if $K=\CO(1)$ is the wave number of a pattern for the rWRIBL, 
then $k=\frac{2\pi}{\lambdahat} K$ is the wave number in the dimensional 
Navier--Stokes system. Thus, if for instance $\hhat=\CO(1)$ and $\ahat=\CO(1)$ 
are fixed such that $\del=2\pi\hhat/\lambdahat=\CO(\eps)$ 
and $\zeta=2\pi\ahat/\lambdahat=\CO(\eps)$ are small due 
to $\lambdahat=\CO(\eps^{-1})$, then $k=\CO(\eps)$.  
However, even in this case, as already 
said at the start of \S\ref{onum}, in applications we always have finite 
$\eps$. For instance, in Fig.~\ref{rc-fig-45} (b)--(d) we 
find $k=\pi/10, k=\pi/8, k=3\pi/20$ [mm$^{-1}$] as the basic 
wave numbers 
(the smallest possible wave number over a domain of length 
80$\,$mm being $\pi/40$ [mm$^{-1}$]).

Over flat bottom, short wave instabilities are only known for
very small inclination angles, see \cite[Section 2.3]{cd02}. In
particular,  calculating the eigenvalues of the
linearization of the rWRIBL \eqref{FTR}, \eqref{QTR} around the Nusselt
solution for the above parameters but $\zeta = 0$ 
by a Fourier ansatz we find no short wave instability in case
of a flat bottom. 

\begin{figure}[!ht]
  \begin{center}
    \begin{tabular}[h]{m{0.0cm}m{6cm}}
      \begin{sideways}
        \footnotesize $u$ [mm/s]
      \end{sideways}
      & \ig{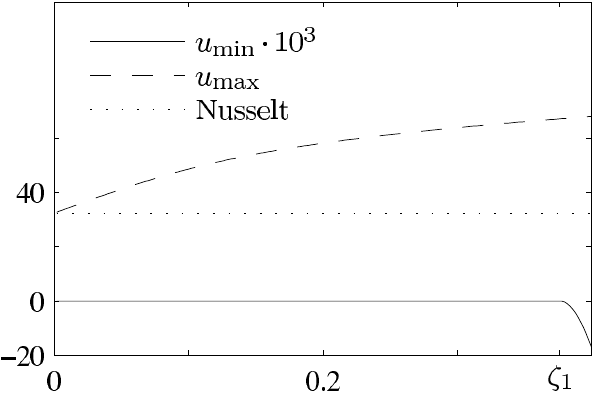} \vspace{-0.5cm} \\
      \small (a) & \centering \hspace{5mm}\footnotesize $\zeta$
    \end{tabular}\hspace{0.4cm}
    \begin{tabular}[h]{m{0.2cm}m{6cm}}
      \begin{sideways}
        \footnotesize $\zhat$ [mm]
      \end{sideways}
      & \ig{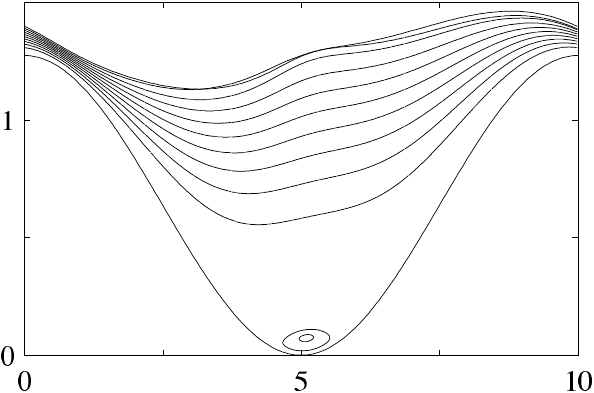} \vspace{-0.5cm} \\
      \small (b) & \centering \footnotesize $\xhat$ [mm]
    \end{tabular}
    \caption{Stationary solutions for parameter set C and $\rey =
      4.2$.  (a) Minimal and maximal downstream velocity
      $u_\text{min}, u_\text{max}$ of stationary solutions depending
      on $\zeta$. For $\zeta>\zeta_1 \approx 0.38 $ eddies occur.  (b)
      Free surface and reconstructed streamlines for $\zeta=0.4$.}
    \label{rc-fig-67}
  \end{center}
\end{figure}

In Fig.~\ref{rc-fig-67}, for fixed $\rey$ and varying $\zeta$ we plot
the minimal and maximal downstream velocities of some stationary
solutions used in Fig.~\ref{rc-fig-23}, which shows that these are
continuations of the Nusselt solution.  For $\zeta> \zeta_1$ the
minimal velocity $u_\text{min}$ becomes negative which is an easy
diagnostic for the existence of eddies. In particular, from
$\zeta_0<\zeta_1$ we find that the short wave instability sets in
before the appearance of eddies, which shows that the short wave
instability is an effect of the wavy bottom on a Nusselt like laminar
solution.  Figure \ref{rc-fig-67} (b) shows the stationary solution
and reconstructed streamlines for the short wave unstable parameters
$\zeta=0.4, \rey=4.2$.

\begin{figure}[!ht]
  \begin{center}\small
  \begin{tabular}[h]{b{0.1cm}@{}b{7cm}@{}}
    (a) & \ig{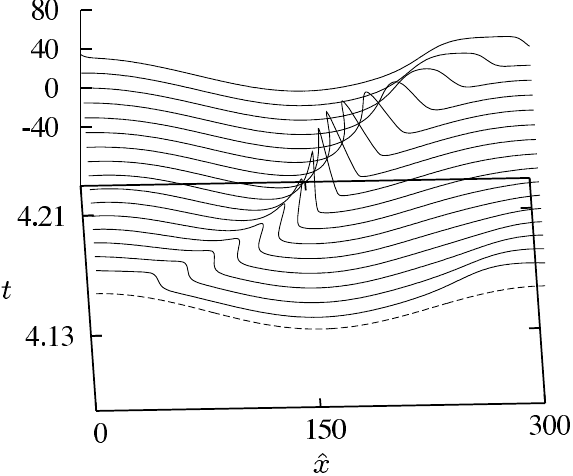}  
  \end{tabular}\hspace{0.1cm}
  \begin{tabular}[h]{b{0.1cm}@{}b{7cm}@{}}
    (b) & \ig{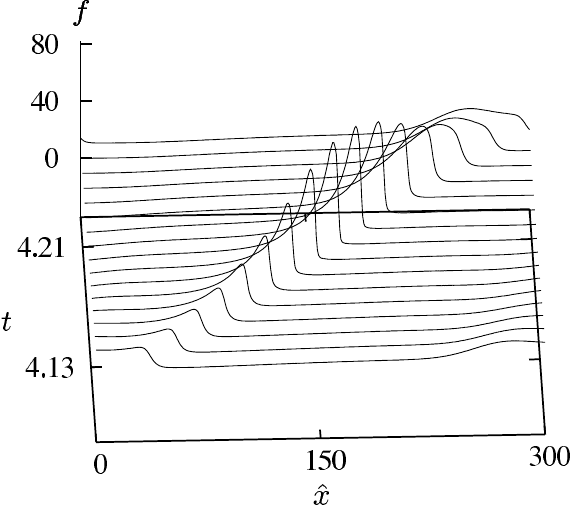}  
  \end{tabular}
    \caption{$\al=90^\circ$, $\lambdahat=\unit{300}{\milli\meter},
      \ahat=\unit{20}{\milli\meter}$, hence $\zeta=0.42$, $\rey=10$,
      $\delta=0.32$; $\ibond=0.003$ and initial data $(F,Q)\equiv
      (1,1)$. (a) Free surface over $\xhat$, dashed line is the bottom
      contour. (b) Film thickness $f$ over $\xhat$.\label{ngf}}
  \end{center}
\end{figure}

Finally, Fig.~\ref{ngf} illustrates a rather strongly unstable
situation where due to a relatively large traveling pulse 
the free surface is not a graph over $\xhat$. This was one of
the motivations to use curvilinear coordinates. Downstream the bottom
maxima where the local inclination angle is larger than $90^\circ$
($\unit{0}{\milli\meter}\ {<}\ \xhat\ {<}\ \unit{150}{\milli\meter}$
in Fig.~\ref{ngf}) we find a bearing-out of the free surface as a
pulse passes.  This overhang is typically rather small since the pulse
is small as it lost mass when it climbed ``uphill''
($\unit{150}{\milli\meter}\ {<}\ \xhat \ {<} \
\unit{300}{\milli\meter}$ in Fig.~\ref{ngf}) to the maximum of the
bottom.  On the other hand, running ``downhill'', the pulse grows and
reaches maximum amplitude around $\xhat\approx
\unit{180}{\milli\meter}$. This yields an overhang (to the left) of
the free surface at the beginning of the ``uphill'' section.

In the literature we did not find data or solutions comparable to
Fig.~\ref{ngf}, or to the short wave instability explained in Figures
\ref{rc-fig-23} to \ref{rc-fig-67}.  Thus we think it will be
interesting to study either experimentally or by full Navier--Stokes
numerics the accuracy of these predictions.

\subsection{Conclusions}\label{ssec}
Using a Galerkin method with only one ansatz and test function we derived the
WRIBL \reff{FT2}, \reff{QT} for film flow over wavy bottom, which in the limit
of flat bottom equals the (one-dimensional version of the) WRIBL derived in
\cite{Scheid_06}. In a second step we regularized the WRIBL to the 
rWRIBL  \reff{FTR}, \reff{QTR}. 
Numerical simulations of the rWRIBL show very good agreement with available
data from experiment and full scale Navier--Stokes numerics.
Finally, our rWRIBL predicts two qualitatively new 
phenomena, namely a short wave instability of Nusselt like solutions 
(without eddies) at non-small inclination angles and at still rather small 
$\zeta$, and solutions where the free surface is not a graph over the 
(Cartesian) downstream coordinate. It remains to be seen whether these 
predictions can be verified experimentally or by full Navier--Stokes numerics.

\begin{appendix}
\section{Curvilinear coordinates}\label{sec_CC}
In order to expand the non-dimensional curvature $K$ and the local
inclination angle $\theta$ in powers of $\zeta$ we first scale the Cartesian
coordinate $\xhat$ and the bottom profile $\bhat$ by
\begin{equation*}
    \Xhat = \frac{2\pi}{\lambdahat}\xhat, \quad \Bhat(\Xhat) = 
    \frac{1}{\ahat}\bhat\left( \frac{\lambdahat}{2\pi}\Xhat \right).
\end{equation*}
This implies $\dxhat\bhat(\xhat) = \zeta \dXhat\Bhat\left(
\frac{2\pi}{\lambdahat}\xhat \right)$ and $ \dxhat[2]\bhat(\xhat) =
\frac{4\pi^2\ahat}{\lambdahat^2}\dXhat[2]\Bhat \left(
\frac{2\pi}{\lambdahat}\xhat \right).  $ The relation between 
$\Xhat$ and $X$ is 
\begin{align*}
X
  & =\frac{2\pi}{\lambdahat}x = \frac{2\pi}{\lambdahat}\int_0^{\xhat}
  \sqrt{1+\left( \dxhat\bhat(\xhat) \right)^2}d\xhat \\
  & = \int_0^{\Xhat}\sqrt{1+\zeta^2\left( \dXhat\Bhat(\Xhat) \right)^2}d\Xhat \\
  & =  \Xhat + \frac{1}{2}\zeta^2\int_0^{\Xhat}\left( \dXhat\Bhat(\Xhat) \right)^2
  d\Xhat + \ord(\zeta^4),
\end{align*}
thus
\mbox{$ \Xhat(X) = X -
  \frac{1}{2}\zeta^2\int_0^{X} \left( \dXhat\Bhat(\Xhat) \right)^2 d\Xhat +
  \ord(\zeta^4),$} 
and therefore $K(X)$ reads (cf.~\eqref{kappa})
  \begin{align}
    K(X) & = \frac{\lambdahat^2}{4\pi^2 \ahat}\kappa\left(
      \frac{\lambdahat}{2\pi}\Xhat(X) \right)\notag\\ 
&= - \frac{\lambdahat^2}{4\pi^2
      \ahat}\frac{\dxhat[2]\bhat\left( \frac{\lambdahat}{2\pi} \Xhat(X)
      \right)}{\left[ 1 + \left( \dxhat\bhat \left( \frac{\lambdahat}{2\pi}
            \Xhat(X) \right) \right)^2 \right]^{\frac{3}{2}}}
    = - \frac{\dXhat[2]\Bhat(\Xhat(X))}
    {\left[ 1 + \zeta^2\left( \dXhat\Bhat(\Xhat(X)) \right)^2 \right]^{\frac{3}{2}}} \nonumber\\
    & = - \dXhat[2]\Bhat(X) + \frac{1}{2}\zeta^2\left( 3\dXhat[2] \Bhat(X)
      (\dXhat\Bhat(X))^2 +
      \dXhat[3]\Bhat(X)\int_0^X(\dXhat\Bhat(\Xhat))^2d\Xhat
    \right) + \ord(\zeta^4) \nonumber\\
    & =: K_0(X) + \zeta^2 K_2(X) + \ord(\zeta^4). \label{K}
  \end{align}
  For the local inclination angle $\theta$ we get
  \begin{align}
    \theta(X) & = \arctan\left(\dxhat\bhat\left(
        \frac{\lambdahat}{2\pi}\Xhat(X) \right) \right) = \arctan(\zeta
    \dXhat\Bhat(\Xhat(X))) = \zeta \dXhat\Bhat(X)
    + \ord(\zeta^3) \nonumber\\
    & =: \zeta \theta_1(X) + \ord(\zeta^3). \label{theta}
  \end{align}

\section{Second-order velocity profile}\label{sec_u2}
Calculating the second-order component of the downstream velocity $U = U_0
+ \eps U_1 + \eps^2 U_2 + \ord(\eps^3)$ by exactly the same approach as in
Section \ref{sec_prof} yields \small
\begin{align}
    \eps^2 U_2 = \ & \delta^2\rey^2\dX[2]F\left( {-}\frac{27}{4480}FZ^8 {+}
      \frac{27}{560}F^2 Z^7 {-} \frac{3}{20}F^3 Z^6 {+} \frac{9}{40}F^4 Z^5
      {+} \frac{3}{8}F^5 Z^4 {-} \frac{21}{10}F^6 Z^3
      {+} \frac{30}{7}F^8 Z \right) \nonumber\\
    & {+} \delta^2\rey^2(\dX F)^2\left( {-}\frac{27}{4480}Z^8 {+}
      \frac{27}{560}FZ^7
      {-} \frac{21}{80}F^2Z^6 {+} \frac{9}{10}F^3Z^5 {+} \frac{15}{8}F^4Z^4 {-} \frac{63}{5}F^5Z^3 {+} \frac{948}{35}F^7Z \right) \nonumber\\
    & {+} \delta^2 \rey(\ibond\dX[4]F {-} \cot\alpha \ \dX[2]F)\left(
      \frac{1}{40}Z^6
      {-} \frac{3}{20}FZ^5 {+} \frac{3}{4}F^2Z^4 {-} 2F^3Z^3 {+} \frac{18}{5}F^5Z \right) \nonumber\\
    & {-} \delta\zeta\rey(\ibond \dX[2]K_0 {+} \cot\alpha \ \dX\theta_1)\left(
      \frac{1}{40}Z^6
      {-} \frac{3}{20}FZ^5 {+} \frac{3}{8}F^2Z^4 {-} \frac{1}{2}F^3Z^3 {+} \frac{3}{5}F^5Z \right) \nonumber\\
    & {+} \delta\rey\left(3\delta\ibond(\dX[2]F)^2 {-} 2\delta\cot\alpha(\dX
      F)^2 {+} 5\delta\ibond\dX[3] F \dX F {-} \zeta(\ibond\dX K_0 {+}
      \cot\alpha \ \theta_1) \dX F\right)
    \cdot \nonumber\\
    & \cdot \left( \frac{3}{4}FZ^4 {-} 3F^2Z^3 {+} 6F^4Z \right)
    {+} \delta^2\rey\ibond(\dX F)^2\dX[2]F\left( \frac{9}{2}Z^4 {-} 18FZ^3 {+} 36F^3Z \right) \nonumber \\
    & {+} \delta\zeta K_0\left( \frac{1}{2}Z^3 {-} \frac{3}{2}FZ^2 {+} 3 F^2Z
    \right) {+} \delta\zeta \dX \theta_1\left( {-}\frac{1}{2}Z^3 {+}
      \frac{3}{2}FZ^2 {-} \frac{3}{2} F^2Z \right)
    \nonumber \\
    & {+} \delta^2\dX[2]F\left( {-}Z^3 {-} \frac{3}{2}FZ^2 {+}
      \frac{15}{2}F^2Z \right) {+} \delta\zeta\theta_1 \dX F \left(
      \frac{3}{2}Z^2 {-} 3FZ \right)\nonumber\\
& {+} \delta^2(\dX F)^2\left({-}\frac{3}{2}Z^2 {+} 15FZ \right) 
    {+} \zeta^2\theta_1^2\left( \frac{3}{4}Z^2 {-} \frac{3}{2}FZ
    \right). \label{U2}
  \end{align}
  \normalsize 
If the flow rate $Q$ is assumed to be enslaved by the film
  thickness $F$, then integration of \eqref{U2} along $Z \in [0,F]$ gives the
  second-order component of $Q = Q_0 + \eps Q_1 + \eps^2 Q_2 + \ord(\eps^3)$, 
see \reff{q2}. 
  In order to achieve an accurate velocity profile $\Uprof_2$ depending on
  both $F$ and $Q$ we again treat $Q$ as independent $\ord(1)$-quantity. This
  profile is not needed for the Galerkin method but only for reconstructing
  flow fields. Therefore we restrict our calculations to the practically
  relevant case of stationary flow over a sinusoidal bottom $\bhat(\xhat) =
  \ahat \cos\left( \frac{2\pi}{\lambdahat}\xhat \right)$. This implies according to
  \eqref{K} and \eqref{theta}
  \begin{equation*}
    K(X) = \cos X + \ord(\zeta^2),\quad \theta(X) = -\zeta\sin X + \ord(\zeta^3).
  \end{equation*}
  In case of the first-order profile \eqref{3_40} the correction of the
  parabolic profile turned out to be a self-similar polynomial with a
  coefficient depending on $Q$ and $\dX Q$ but not on $F$ or its spatial
  derivatives. The basic assumption now is that this is also true for the
  second-order correction $\Uprof_2$. Therefore all spatial derivatives of $F$
  emanate from $\dX Q$. As we consider here only stationary solutions the
  evolution equation for $F$ \eqref{FT} gives $\dX Q = 0$. Thus in $U_2$ we
  neglect all terms containing spatial derivatives of $F$. Taking again into
  account that treating $Q$ as independent quantity mixes up $\eps$-orders in
  the expansion of $U$ finally yields\small
  \begin{align*}
    \eps^2\Uprof_2 = \ & \delta\zeta\rey(\cot\alpha + \ibond)\cos X Q^2 \left(
      \frac{1}{40} \left( \frac{Z}{F} \right)^6 {-} \frac{3}{20}\left(
        \frac{Z}{F} \right)^5 + \frac{1}{4}\left( \frac{Z}{F} \right)^4 {-}
      \frac{9}{35}\left( \frac{Z}{F} \right)^2
      + \frac{4}{35}\frac{Z}{F} \right) \nonumber\\
    & + \delta\zeta\cos X Q \left( \left( \frac{Z}{F} \right)^3
      {-}\frac{3}{4}\left( \frac{Z}{F} \right)^2 \right).
  \end{align*}
\normalsize
  Thus the velocity profile used in Section \ref{sec_num} to reconstruct flow
  fields reads\small
  \begin{align}
    \Uprof = \ & \frac{3Q}{F}\left( {-}\frac{1}{2}
      \left(\frac{Z}{F}\right)^2+\frac{Z}{F} \right) + \delta\zeta\cos X Q
    \left( \left( \frac{Z}{F} \right)^3
      {-}\frac{3}{4}\left( \frac{Z}{F} \right)^2 \right)  \nonumber \\
    & + \delta\zeta\rey(\cot\alpha + \ibond)\cos X Q^2 \left( \frac{1}{40}
      \left( \frac{Z}{F} \right)^6 {-} \frac{3}{20}\left( \frac{Z}{F}
      \right)^5 + \frac{1}{4}\left( \frac{Z}{F} \right)^4 {-}
      \frac{9}{35}\left( \frac{Z}{F} \right)^2 + \frac{4}{35}\frac{Z}{F}
    \right).
    \label{Uprof2}
  \end{align}
\normalsize
  The according velocity component $\Wprof$ is given by the continuity
  equation, i.e.
  \begin{equation}
    \Wprof = -\frac{1}{1+\delta\zeta \cos X \ Z}\int_0^Z \dX U dZ.
    \label{Wprof2}
  \end{equation}

\end{appendix}

\noindent
{\bf Acknowledgement:} This work was supported by the DFG under grant Schn
520/6. The authors thank Andreas Wierschem and Vasilis Bontozoglou for
stimulating discussions during early stages of this work. 

\bibliographystyle{alpha}

\small

\bibliography{./literatur.bib}

\newcommand{\etalchar}[1]{$^{#1}$}
\begin{thebibliography}{VMHM05}

\bibitem[AVB06]{avb06}
K.~Argyriadi, M.~Vlachogiannis, and V.~Bontozoglou.
\newblock Experimental study of inclined film flow along periodic corrugations:
  The effect of wall steepness.
\newblock {\em Phys. Fluids}, 18:012102, 2006.

\bibitem[Ben57]{Benjamin_57}
T.~B. Benjamin.
\newblock Wave formation in laminar flow down an inclined plane.
\newblock {\em J. Fluid Mech.}, 2:554--574, 1957.

\bibitem[Ben66]{Benney_66}
D.~J. Benney.
\newblock Long waves on liquid films.
\newblock {\em J. of Mathematics and Physics}, 45:150--155, 1966.

\bibitem[CD02]{cd02}
H.-C. Chang and E.~A. Demekhin.
\newblock {\em Complex Wave Dynamics on Thin Films}.
\newblock Elsevier, Amsterdam, 2002.

\bibitem[DO07]{dav07}
L.~A. D\'avalos-Orozco.
\newblock Nonlinear instability of a thin film flowing down a smoothly deformed
  surface.
\newblock {\em Phys. Fluids}, 19:074103, 2007.

\bibitem[HBAW09]{hbaw08}
C.~Heining, V.~Bontozoglou, N.~Aksel, and A.~Wierschem.
\newblock {Nonlinear resonance in viscous films on inclined wavy planes}.
\newblock {\em {Int. J. Multiphase Flow}}, {35}({1}):{78--90}, {2009}.

\bibitem[Lin74]{Lin_74}
S.~P. Lin.
\newblock Finite amplitude side-band stability of a viscous film.
\newblock {\em J. Fluid Mech.}, 63(3):417--429, 1974.

\bibitem[OGN08]{ogn08}
A.~Oron, O.~Gottlieb, and E.~Novbari.
\newblock Numerical analysis of a weighted-residual integral boundary-layer
  model for nonlinear dynamics of falling liquid films.
\newblock {\em European Journal of Mechanics - B/Fluids}, 2008.
\newblock In press.

\bibitem[OH08]{oh08}
A.~Oron and C.~Heining.
\newblock Weighted-residual integral boundary-layer model for the nonlinear
  dynamics of thin liquid films falling on an undulating vertical wall.
\newblock {\em Phys. Fluids}, 20:082102, 2008.

\bibitem[Oos99]{Ooshida99}
T.~Ooshida.
\newblock {Surface equation of falling film flows with moderate Reynolds number
  and large but finite Weber number}.
\newblock {\em {Phys. Fluids}}, {11}({11}):{3247--3269}, {1999}.

\bibitem[PMP83]{Pumir_83}
A.~Pumir, P.~Manneville, and Y.~Pomeau.
\newblock On solitary waves running down an inclined plane.
\newblock {\em J. Fluid Mech.}, 135:27--50, 1983.

\bibitem[PN03]{pano03}
C.~D. Park and T.~Nosoko.
\newblock Three-dimensional wave dynamics on a falling film and associated mass
  transfer.
\newblock {\em AIChE Journal}, 49(11):2715--2727, 2003.

\bibitem[Poz88]{poz88}
C.~Pozrikidis.
\newblock The flow of a liquid film along a periodic wall.
\newblock {\em J. Fluid. Mech.}, 188:275--300, 1988.

\bibitem[RQM98]{Ruyer_98}
C.~Ruyer-Quil and P.~Manneville.
\newblock Modeling film flows down inclined planes.
\newblock {\em Eur. Phys. J. B}, 6:277--292, 1998.

\bibitem[RQM00]{Ruyer_00}
C.~Ruyer-Quil and P.~Manneville.
\newblock Improved modeling of flows down inclined planes.
\newblock {\em Eur. Phys. J. B}, 15:357--369, 2000.

\bibitem[SAB94]{Salamon94}
T.~R. Salamon, R.~C. Armstrong, and R.~A. Brown.
\newblock {Traveling waves on vertical films -- numerical analysis using the
  finite-element method}.
\newblock {\em {Phys. Fluids}}, {6}({6}):{2202--2220}, {1994}.

\bibitem[Shk67]{Shkadov_67}
V.~Y. Shkadov.
\newblock Wave conditions in the flow of a thin layer of a viscous liquid under
  the action of gravity.
\newblock {\em Izv. Akad. Nauk. SSSR, Mekh. Zhidk. Gaza}, 2:43--51, 1967.

\bibitem[SRQM06]{Scheid_06}
B.~Scheid, C.~Ruyer-Quil, and P.~Manneville.
\newblock Wave patterns in film flows: modelling and three-dimensional waves.
\newblock {\em J. Fluid Mech.}, 562:183--222, 2006.

\bibitem[Tri98]{Trifonov_98}
Y.~Y. Trifonov.
\newblock Viscous liquid film flows over a periodic surface.
\newblock {\em Int. J. Multiphase Flow}, 24(7):1139--1161, 1998.

\bibitem[Tri04]{trif04}
Y.~Y. Trifonov.
\newblock Viscous film flow down corrugated surfaces.
\newblock {\em J. Appl. Mech. and Techn. Phys.}, 45(3):389--400, 2004.

\bibitem[Tri07a]{trif07a}
Y.~Y. Trifonov.
\newblock {S}tability and nonlinear wavy regimes in downward film flows on a
  corrugated surface.
\newblock {\em J. Appl. Mech. and Techn. Phys.}, 48(1):91--100, 2007.

\bibitem[Tri07b]{trif07}
Y.~Y. Trifonov.
\newblock Stability of a viscous liquid film flowing down a periodic surface.
\newblock {\em Intern. J. of Multiphase Flow}, 33:1186--1204, 2007.

\bibitem[VB02]{vb02}
M.~Vlachogiannis and V.~Bontozoglou.
\newblock Experiments on laminar film flow along a periodic wall.
\newblock {\em J. Fluid Mech.}, 457:133--156, 2002.

\bibitem[VMHM05]{vmhm05}
O.~K. Valluri, G.~Matar, F.~Hewitt, and M.~A. Mendes.
\newblock Thin film flow over structured packings at moderate {R}eynolds
  numbers.
\newblock {\em Chem. Eng. Sci.}, 60:1965--1975, 2005.

\bibitem[WA03]{wa03}
A.~Wierschem and N.~Aksel.
\newblock Instability of a liquid film flowing down an inclined wavy plane.
\newblock {\em Phys. D}, 186(3--4):221--237, 2003.

\bibitem[WBH{\etalchar{+}}08]{wbhua08}
A.~Wierschem, V.~Bontozoglou, C.~Heining, H.~Uecker, and N.~Aksel.
\newblock Linear resonance in viscous films on inclined wavy planes.
\newblock {\em Int. J. Multiphase Flow}, 34:580--590, 2008.

\bibitem[WLA05]{Acta_05}
A.~Wierschem, C.~Lepski, and N.~Aksel.
\newblock Effect of long undulated bottoms on thin gravity-driven films.
\newblock {\em Acta Mech.}, 179:41--66, 2005.

\bibitem[WSA03]{wsa03}
A.~Wierschem, M.~Scholle, and N.~Aksel.
\newblock Vortices in film flow over strongly undulated bottom profiles at low
  {R}eynolds numbers.
\newblock {\em Phys. Fluids}, 15(2):426--435, 2003.

\bibitem[Yih63]{Yih_63}
C.~Yih.
\newblock Stability of liquid flow down an inclined plane.
\newblock {\em Phys. Fluids}, 6(3):321--334, 1963.

\end{thebibliography}

\end{document}